\documentclass[11pt]{article}
\input epsf.tex
\usepackage{amssymb}
\usepackage{latexsym}
\font\tensym=msbm10
\font\sevensym=msbm7
\font\fivesym=msbm5

\newfam\ssymfam
\textfont\ssymfam=\tensym
\scriptfont\ssymfam=\sevensym
\scriptscriptfont\ssymfam=\fivesym

\catcode`\@=11
\renewcommand\subsection{\@startsection{subsection}{2}{\z@}%
                                     {-3.25ex\@plus -1ex \@minus -.2ex}%
                                     {-0.01 mm}
                                     {\normalfont\large\bfseries}}

\catcode`\@=11
\renewcommand\subsubsection{\@startsection{subsubsection}{2}{\z@}%
                                     {-3.25ex\@plus -1ex \@minus -.2ex}%
                                     {-0.01 mm}
                                     {\normalfont\bfseries}}

\newtheorem{example}{Example}
\newtheorem{theorem}[example]{Theorem}
\newtheorem{corollary}[example]{Corollary}
\newtheorem{proposition}[example]{Proposition}
\newtheorem{lemma}[example]{Lemma}
\newtheorem{conjecture}[example]{Conjecture}

\def\resp{{\em resp.$\ $}}
\def\proof{\medskip\noindent {\it Proof --- \ }}
\def\finex{\hfill $\Diamond$}

\def\cqfd{\hfill $\Box$ \bigskip}
\def\adots{\mathinner{\mkern2mu\raise1pt\hbox{.}
\mkern3mu\raise4pt\hbox{.}\mkern1mu\raise7pt\hbox{.}}}
\def\<{\lngle\,}
\def\>{\,\rangle}

\def\cf{{\it cf.$\ $}}

\def\hom{{\rm Hom}}
\def\ext{{\rm Ext}}

\def\ie{{\it i.e. }}

\def\tab{{\rm Tab\, }}

\def\inv{{\rm inv}}
\def\SG{\mathfrak S}

\def\a{\alpha}

\def\b{\beta}
\def\N{{\bf N}}
\def\Z{{\bf Z}}
\def\C{{\bf C}}

\def\Q{{\bf Q}\, }

\def\F{{\cal F}}

\def\B{{\cal B}}

\def\V{{\mathfrak V}}

\def\I{{\cal I}}
\def\J{{\cal J}}
\def\SS{{\cal S}}

\def\mod{{\rm \ mod\ }}

\def\wt{{\rm wt}}
\def\g{\mathfrak g}

\def\gl{\mathfrak gl}
\def\Sl{\mathfrak sl}
\def\Rg{\mathfrak R}
\def\nn{\mathfrak n}
\def\bb{\mathfrak b}
\def\slchap{\widehat{\mathfrak sl}}

\def\H{\widehat H}

\def\A{{\cal A}}
\def\AA{{\bf A}}

\def\L{{\cal L}}

\def\bar{\overline}

\def\<{\langle}
\def\>{\rangle}

\def\deg{{\rm deg}}
\def\BB{{\cal B}}
\def\CC{{\cal C}}
\def\m{{\bf m}}
\def\n{{\bf n}}
\def\p{{\bf p}}
\def\q{{\bf q}}
\def\a{{\bf a}}
\def\b{{\bf b}}
\def\s{{\bf s}}
\def\t{{\bf t}}
\def\M{{\cal M}}
\def\RR{{\cal R}}

\def\le{\leqslant}
\def\ge{\geqslant}

\newdimen\Squaresize \Squaresize=14pt
\newdimen\Thickness \Thickness=0.5pt
 
\def\Square#1{\hbox{\vrule width \Thickness
   \vbox to \Squaresize{\hrule height \Thickness\vss
      \hbox to \Squaresize{\hss#1\hss}
   \vss\hrule height\Thickness}
\unskip\vrule width \Thickness}
\kern-\Thickness}
 
\def\Vsquare#1{\vbox{\Square{$#1$}}\kern-\Thickness}

\def\young#1{
\vbox{\smallskip\offinterlineskip
\halign{&\Vsquare{##}\cr #1}}}
 
\title{\bf Induced representations of affine Hecke algebras
and canonical bases of quantum groups}
\author{\rm Bernard {\sc Leclerc},
\ Maxim {\sc Nazarov}
\rm \ and Jean-Yves {\sc Thibon}
}

\date{}

\begin{document}
\maketitle

\vskip 1cm

\begin{abstract}
A criterion of irreducibility for induction products
of evaluation modules of type~$A$ affine Hecke algebras
is given.
It is derived from multiplicative properties of
the cano\-ni\-cal basis of a quantum deformation of the 
Bernstein-Zelevinsky ring.
\end{abstract}

\vskip 0.6cm

\section{Introduction} \label{SECT1}
Let $GL_m=GL(m,F)$ be the general linear group over
a non-Archimedean local field $F$, and $B_m$ its
subgroup of upper triangular matrices. 
To each $m$-tuple of complex para\-me\-ters $(s_1,\ldots ,s_m)$ is
attached a representation $I(s_1,\ldots ,s_m)$ of $GL_m$ 
in the unramified principal series,
which is defined as 
the right regular representation of $GL_m$ on the space of locally
constant functions $f : GL_m \longrightarrow \C$ satisfying
\[
f(b\,g) = \chi(b)\,f(g), 
\qquad (b\in B_m,\,g\in GL_m).
\] 
Here, $\chi$ is the character of $B_m$ given on $b=(b_{ij})$ by
\[
\chi(b)=\delta(b)^{1\over 2}\,\prod_{i=1}^m |b_{ii}|_F^{s_i}\,,
\]
and 
\[
\delta(b) = \prod_{1\le i<j \le m} {|b_{ii}|_F\over|b_{jj}|_F}
= \prod_{i=1}^m |b_{ii}|_F^{m+1 - 2i}
\]
is the modulus of $B_m$.
We denote by $\CC_m$ the category of smooth representations of finite
length of $GL_m$ whose composition factors are subquotients of 
$I(s_1,\ldots , s_m)$ for some choice of $(s_1,\ldots ,s_m)\in\C^m$.

Let ${\mathfrak o}$ be the ring of integers of $F$, ${\mathfrak p}$ its maximal
ideal, $k={\mathfrak o}/{\mathfrak p}$ its residue field of cardinality $q$.
Let $I_m$ be the standard Iwahori subgroup of $GL_m$ consisting of those
$g=(g_{ij}) \in GL(m,{\mathfrak o})$ for which 
$g_{ij}\in {\mathfrak p}$ whenever $i>j$. 
The space $\H_m$ of compactly supported complex-valued functions on $GL_m$ 
which are bi-invariant with respect to $I_m$ is an algebra under
convolution called the Iwahori-Hecke algebra of $GL_m$.

By a theorem of Bernstein, Borel, Casselman and Matsumoto
(see e.g. \cite{Wa}), the category $\CC_m$
is equivalent to the category $\CC(\H_m)$ of 
finite-dimensional complex representations of $\H_m$.
Moreover, one has natural induction functors 
$(M_1,M_2) \mapsto M_1\odot M_2$ from 
$\CC_{m_1} \times \CC_{m_2}$ to $\CC_{m_1+m_2}$ and from
$\CC(\H_{m_1}) \times \CC(\H_{m_2})$ to $\CC(\H_{m_1+m_2})$
which correspond to each other via the equivalences of categories.
In this paper we want to study under which conditions induction
products $L_1 \odot L_2$ of irreducible objects of these categories 
are irreducible.

Following Bernstein and Zelevinsky, one introduces 
\[\RR = \bigoplus_{m\ge 0} \RR_m,\]
where for $m\ge 1$, $\RR_m$ denotes the complexified Grothendieck group 
of $\CC_m$ (or $\CC(\H_m)$), and $\RR_0 := \C$. 
The induction and restriction functors make $\RR$ into a Hopf algebra
which was explicitely described by Zelevinsky \cite{Zel80}.
For $s\in\C$, let $\RR_m(s)$ be the subgroup of $\RR_m$ generated by
the classes of the subquotients of all representations $I(s_1,\ldots ,s_m)$
for which $s_1,\ldots ,s_m \in s+\Z$, and let 
\[\RR(s) = \bigoplus_{m\ge 0} \RR_m(s).\]
Note that, since $|x|_F^s=q^{-s\, {\rm val}_F(x)}$ depends only on 
$s$ modulo $(2\pi i/{\rm log}q)\Z$, the group $\RR(s)$ depends only on
$s$ modulo $\Omega = \Z\oplus (2\pi i/{\rm log}q)\Z$.
Now (\cite{Zel80}, 8.7), $\RR(s)$ is a subalgebra of $\RR$, the $\RR(s)$ are all isomorphic to each
other, and $\RR = \bigotimes \RR(s)$ where $s$ runs over the elliptic curve
$\C/\Omega$.
So it is enough to describe $R := \RR(0)$.

It turns out that $R$ is isomorphic, as a Hopf algebra, to the algebra
of polynomials in the coordinate functions of the group $N_\infty$ of upper triangular
unipotent $\Z\times\Z$-matrices with finitely many non-zero entries
off the main diagonal. As such, $R$ has a natural quantum deformation
$R_v$ (in the sense of Drinfeld and Jimbo).
By results of Kashiwara and Lusztig, $R_v$ is endowed with
a canonical basis $\BB_v$ which specializes when $v\mapsto 1$ to a
canonical basis $\BB$ of $R$. Our approach to the representation theory
of $GL_m$ and $\H_m$ is based on the following crucial fact:
\begin{quote}
{\em The canonical basis $\BB$ coincides with the basis of $R$ consisting of
the classes of irreducible representations.}
\end{quote}
This theorem follows by comparing the $p$-adic analogue of the Kazhdan-Lusztig
formula, conjectured by Zelevinsky \cite{Zel81} and proved by Ginzburg \cite{CG}, 
with the geometrical description by Lusztig of the canonical basis
\cite{Lu90}.
The dual version of this theorem is proved in \cite{Ar}.
 
Here $\BB_v$ is a dual canonical basis in the sense of Lusztig, or an
upper global basis in the sense of Kashiwara.
A distinguished subset of $\BB_v$ consists of the quantum flag minors
\cite{BZ,LZ}. The quantum flag minors of degree $m$ correspond 
to a special class of irreducible
$\H_m$-modules, called the evaluation modules. 
They are obtained by lifting the simple 
modules of the
finite-dimensional Hecke algebra $H_m$ via the evaluation maps
(see below \S \ref{SECT13}). 
In the equivalence of categories $\CC(\H_m) \simeq \CC_m$, these
modules correspond in turn to a large class of irreducible
representations of $GL_m$,
namely those parametrized, up to a shift by an arbitrary $z\in\C^*$,
by Zelevinsky's
multi-segments of the form
\[
\m = \sum_{k=1}^r \, [i+k-1, j_k],
\]
where $i,j_1,\ldots , j_k$ are integers, and  
$i< j_1 < j_2 < \cdots <j_r$.
 
The theorem above implies that multiplying vectors of $\BB$ is the same
as taking induction products of irreducible modules in $R$, and therefore
we are led to the problem of understanding which products of 
elements of $\BB$ belong to $\BB$.
The multiplicative properties of the dual canonical basis have been 
studied by Berenstein and Zelevinsky \cite{BZ}. 
They conjectured that the product of
two vectors of $\BB_v$ 
belongs to $\BB_v$ up to a power of $v$, if and
only if these vectors commute up to a power of~$v$.
We will prove this conjecture in the special case of the 
quantum flag minors, for which we have an explicit criterion
of $v$-commutativity \cite{LZ}. This will give a criterion for the
irreducibility of induction products of evaluation modules.
Namely, for $z\in\C^*$ and $\alpha$ a partition of $m$, 
let $S_\alpha(z)$ denote the $\H_m$-module obtained by evaluation
at $z$ of the simple $H_m$-module $S_\alpha$ attached to a partition
$\alpha$ of $m\ge 1$.
Similarly, let $S_\beta(w)$ be an evaluation module of $\H_n(z)$ for
some $n\ge 1$.
Since the induction product $S_\alpha(z) \odot S_\beta(w)$
is simple if $z/w \not\in q^\Z$ \cite{Zel80},
we can assume that $z/w \in q^{\Z}$.
Our main result is 
\begin{theorem} \label{MAIN1} Let $z/w = q^c$ for some $c\in\Z$.
Associate to the partitions $\alpha = (\alpha_1,\ldots ,\alpha_r)$
and $\beta = (\beta_1, \ldots , \beta_s)$ the following subsets of $\Z$ :
\[
\I = \Z_{\le c-r} \cup \{c-r+1+\alpha_r,\ldots ,c+\alpha_1\}, \qquad
\J = \Z_{\le -s} \cup \{-s+1+\beta_s,\ldots ,\beta_1\}.
\]
{\rm (i)} Suppose $c > 0$. Then, the product $S_\alpha(z) \odot S_\beta(w)$ is
not simple if and only if there exist $i,j,k \in \Z$ such that  
$i,k \in \I\setminus \J$,  $j\in \J\setminus \I$ and $i<j<k$.

\smallskip\noindent
{\rm (ii)} Suppose $c < 0$. Then, the product $S_\alpha(z) \odot S_\beta(w)$ is
not simple if and only if there exist $i,j,k \in \Z$ such that  
$i,k \in \J\setminus \I$,  $j\in \I\setminus \J$ and $i<j<k$.

\smallskip\noindent
{\rm (iii)} Suppose $c = 0$.
Then, the product $S_\alpha(z) \odot S_\beta(w)$ is
not simple if and only if there exist $i,j,k,l \in \Z$ such that
$i,k \in \I\setminus \J$,  $j,l\in \J\setminus \I$  and either
$i<j<k<l$ or $j<i<l<k$.
\end{theorem}
Note that $S_\alpha(z) \odot S_\beta(w)$ is simple if and only if
$S_\beta(w) \odot S_\alpha(z)$ is simple. Hence, the statements
(i) and (ii) of our Theorem \ref{MAIN1} are in fact equivalent. 
Using an argument of Kitanine, Maillet and Terras \cite{KMT, MT} (see also
\cite{NT2}), 
we deduce from Theorem~\ref{MAIN1}
an irreducibility criterion for products of any 
number of evaluation modules,
by translating to representations of affine quantum groups.
\begin{theorem} \label{MAIN2}
Let $z_1,\ldots ,z_r$ be non-zero complex numbers and 
$\alpha^{(1)},\ldots ,\alpha^{(r)}$ be partitions of
some positive integers.
The product 
$S_{\alpha^{(1)}}(z_1) \odot \cdots \odot S_{\alpha^{(r)}}(z_r)$
is simple if and only if the products
$S_{\alpha^{(k)}}(z_k) \odot S_{\alpha^{(l)}}(z_l)$
are simple for all $1\le k < l \le r$.
\end{theorem}

For partitions $\alpha^{(1)},\ldots ,\alpha^{(r)}$ of certain special
types, the irreducibility criterion for the induced module
$S_{\alpha^{(1)}}(z_1) \odot \cdots \odot S_{\alpha^{(r)}}(z_r)$
has been known. When $\alpha^{(k)}=(1)$ for every $k=1,\ldots,r$
the induced module belongs to the principal
series, and  the irreducibility criterion  has
been known for a long time (in \cite{Kato} 
this criterion was given for the affine Hecke algebras
corresponding to arbitrary root systems).
When each of the partitions $\alpha^{(1)},\ldots ,\alpha^{(r)}$
consists of one part only, 
the induced module belongs to the generalized principal
series, and  the irreducibility criterion
was given by Zelevinsky \cite{Zel80} in terms of segments.

We note that Reineke has also proven a partial case 
of the Berenstein-Zelevinsky conjecture, namely when one of the two
vectors is a ``small'' quantum minor \cite{Rei}.

Our paper is organized as follows. 
In Section~\ref{SECT2} we recall the classification of the irreducible
finite-dimensional representations of the affine Hecke algebras, 
the description of the Bernstein-Zelevinsky algebra $R$, and the
formula for the composition multiplicities of the standard induced
modules.
In Section~\ref{SECT3}, we review the definitions of the quantum
algebras $U^+_v=U_v^+(\Sl_\infty)$ and $R_v$ and of their canonical bases.
In Section~\ref{SECT4}, we describe a simple algorithm for calculating
the dual canonical basis of $U_v^+$.
In Section~\ref{SECT5}, we recall the Berenstein-Zelevinsky conjecture
about multiplicative properties of this basis, as well as the
criterion of quasi-commutativity of \cite{LZ} for quantum flag
minors.
In Section~\ref{SECT6}, we state our results about the 
Berenstein-Zelevinsky conjecture for flag minors, and we prove
most of them in Section~\ref{SECT7}. 
Finally in Section~\ref{SECT8}, we return to evaluation modules
of affine Hecke algebras and derive our main results.


\section{Representations of affine Hecke algebras} \label{SECT2}

\subsection{}   
By a theorem of Bernstein, the Iwahori-Hecke algebra $\H_m=\H_m(q)$ 
of $GL_m$ has the following presentation.
It is the associative $\C$-algebra with invertible generators
$y_1,\ldots,y_m$ and $T_1,\ldots,T_{m-1}$  subject to the relations
\[
\begin{array}{ll}
T_iT_{i+1}T_i=T_{i+1}T_iT_{i+1}, \quad & 1\le i\le m-2,\\[1mm]
T_iT_j=T_jT_i,&\vert i-j\vert>1,\\[1mm]
(T_i-q)(T_i+1)=0,& 1\le i\le m-1,\\[1mm]
y_iy_j=y_jy_i, & 1\le i,j\le m,\\[1mm]
y_jT_i=T_iy_j,& j\not= i,i+1,\\[1mm]
T_iy_iT_i=q\,y_{i+1}, & 1\le i\le m-1 .
\end{array}
\]
More generally, one can define for any $t\in\C^*$ an algebra
$\H_m(t)$ by replacing $q$ by $t$ in this presentation.
It is known that the representation theory of
$\H_m(t)$ is the same for all parameters $t$ which are not roots of
unity, so from now on we only assume that $t$ is a complex
number of infinite multiplicative order 
and we write $\H_m$ in place of $\H_m(t)$.
In particular, our proofs of Theorem~\ref{MAIN1} and
Theorem~\ref{MAIN2} will be valid for any such parameter~$t$.
Therefore, when we refer below to Theorem~\ref{MAIN1} and 
Theorem~\ref{MAIN2}, it should be understood that the 
parameter $q$ of Section~\ref{SECT1}
is replaced by such a generic complex number $t$.

\subsection{}\label{SECT12}
The classification of irreducible finite-dimensional representations 
of $\H_m$ was obtained by Ze\-le\-vinsky \cite{Zel80}
(in the case $t=q$; for the general case see for instance
\cite{Rog2}).
Let $\A_m$ denote the commutative subalgebra of $\H_m$ generated by the 
elements $y_i^\pm\ (1\le i \le m)$.
As explained in the introduction, we may restrict ourselves to the
category of finite-dimensional representations in which the
generators $y_i$ of $\A_m$ have all their eigenvalues in $t^\Z$.
Let $\CC_m^\Z$ denote this category.
The parametrization of the simple modules of $\CC_m^\Z$
is in terms of combinatorial objects
called multi-segments.
It is obtained as follows.
For $\mu=(\mu_1,\ldots,\mu_r)$  a composition of $m$,  set
$$D(\mu)=\{\mu_1+\cdots+\mu_k,\ 1\le k\le r-1\}=\{d_1,\ldots,d_{r-1}\}$$
and denote
by $\H_\mu$ the subalgebra of $\H_m$ generated by $\A_m$ and
$\{T_i\ |\ i\not\in D(\mu)\}$.
For $\a=(a_1,\ldots ,a_r)\in\Z^r$, let $\C_{\mu,\a}$ be the 1-dimensional representation
of $\H_\mu$ defined by $T_i\mapsto t$, $i\not\in D(\mu)$, and
$y_{d_{i-1}+1}\mapsto t^{a_i}$,
$i=1,\ldots,r$, where $d_0=0$.
(Observe that because of the defining relations of $\H_m(t)$
the action of the other generators $y_j$ on $\C_{\mu,\a}$ is then given
by $y_{d_{i-1}+k}\mapsto t^{a_i+k-1}$, $k=1,\ldots ,\mu_i$.)
We denote by $M_{\mu,\a}$ the induced $\H_m$-module
\[
M_{\mu,\a} = \H_m \otimes_{\H_\mu} \C_{\mu,\a} \,.
\]
The representations $M_{\mu,\a}$ belong to the so-called generalized principal
series.
In particular, when $\mu = (1,\ldots ,1) = (1^m)$ and $\a = (a_1,\ldots,a_m)$,
we get a principal series representation of $\H_m$.

Let $R_m$ be the complexified Grothendieck group of the category $\CC_m^\Z$.
When $(\nu,\b)$ is a permutation of $(\mu,\a)$, \ie $\nu_i = \mu_{\sigma(i)}$
and $b_i=a_{\sigma(i)}$ for some $\sigma \in \SG_r$, the induced modules
$M_{\mu,\a}$ and $M_{\nu,\b}$ are in general non isomorphic, but
their classes in $R_m$ are equal \cite{Zel80},
that is, the modules $M_{\mu,\a}$ and $M_{\nu,\b}$ have the same composition
factors with identical multiplicities.

The equivalence class of $(\mu,\a)$ considered up to permutations 
is nothing but a {\em multi-segment}.
Here, by a {\em segment} we mean an interval $[i,j]$ in $\Z$, and we call
multi-segment a formal finite unordered sum
$\m = \sum_{i\le j} m_{ij}[i,j]$.
The integer $m_{ij}$ is the multiplicity of the segment $[i,j]$
in $\m$. When $i=j$ we sometimes write $[i]$ instead of $[i,i]$.
We define the degree of $\m$ by 
$\deg(\m)= \sum_{i\le j}m_{ij}(j-i+1)$.
Now, by attaching to $(\mu,\a)$ the multi-segment
$$
\m=\sum_{i=1}^r \,[a_i,a_i+\mu_i-1]\,,
$$
we obviously get a one-to-one correspondence between the classes
of pairs $(\mu,\a)$ as above and the multi-segments $\m$ of
degree $m$.
Thus we can write unambiguously $[M_\m]$ for the class of $M_{\mu,\a}$ in 
$R_m$.
The set of segments and the set of multi-segments will be denoted 
respectively by $\SS$ and $\M$.

Zelevinsky has introduced the following partial order on $\M$ 
\cite{Zel80,Zel81}.
We say that two segments $\s$ and $\s'$ are linked when $\t=\s\cup \s'$ 
is a segment and $\t$ is different from $\s$ and $\s'$.
Let $\m$ and $\n$ be two distinct multi-segments.
We write $\m \rightarrow \n$ if 
the multi-segment $\n$ can be obtained
from $\m$ by replacing a pair $\s, \, \s'$ of linked segments of $\m$ by
the pair $\t = \s\cup \s',\, \t' = \s \cap \s'$,
where $\t'$ is allowed to be empty.
More generally we write $\m \lhd \n$ if there exists a sequence
of multi-segments $\n_1,\ldots ,\n_k$ such that 
$\m \rightarrow \n_1\rightarrow \cdots \rightarrow \n_k\rightarrow\n$.
It is known \cite{Zel80,Zel81}
that $M_\m$ is irreducible if and only if $\m$ is maximal for
$\unlhd$, and if this is not the case, there is a unique simple module 
whose class occurs in the expansion of $[M_\m]$ but does not occur
in any $[M_\n]$ with $\m\lhd \n$.
Let $L_\m$ denote this simple module (if $\m$ is maximal, $L_\m =
M_\m$). Then the $L_\m$ are pairwise
non isomorphic and all simple modules of $\CC_m^\Z$ are of this type
for some multi-segment $\m$ of degree $m$.

\begin{example}{\rm
If $\H_m = \H_m(q)$ is the Iwahori-Hecke algebra of $GL_m$, then
each finite-dimensional representation $M$ of $\H_m$ 
can be identified with the subspace of $I_m$-fixed vectors
in a smooth representation of finite length $V_M$ of $GL_m$.
In this correspondence the principal series representation
$M_{(1^m),\a}$ of $\H_m$ is sent to the principal series 
representation $I(-a_1,\ldots , -a_m)$ of $GL_m$.

For an irreducible representation $M=L_\m$ of $\H_m$, let
$V_\m=V_M$ denote the corresponding irreducible representation
of $GL_m$.
For instance, consider for $i\in\Z$ the multi-segments 
\[
\m(i) = [i\,,\,i+m-1],\qquad \n(i)=[i]+[i+1]+\cdots +[i+m-1].
\]
Then, $L_{\m(i)}$ is the $1$-dimensional $\H_m$-module in which
all $T_i$'s act by multiplication by $t$, and $y_1$ by multiplication 
by $t^i$. On the other hand,  
$L_{\n(i)}$ is the $1$-dimensional $\H_m$-module in which
all $T_i$'s act by multiplication by $-1$, and $y_1$ by multiplication 
by $t^{i+m-1}$.
Correspondingly,
$V_{\m(i)}$ is the 
$1$-dimensional representation attached to the character
\[
g \mapsto |\det g|_F^{-i-{m-1\over2}}
\]
of $GL_m$, and 
$V_{\n(i)} = V_{\m(i)} \otimes {\rm St}$, where ${\rm St}$ denotes the
special (or Steinberg) representation, that is, the unique irreducible
quotient of the right regular representation of $GL_m$ on the space
of locally constant functions on $B_m {\setminus} GL_m$.
}
\finex
\end{example}
\begin{figure}[t]
\begin{center}
\leavevmode
\epsfxsize =10cm
\epsffile{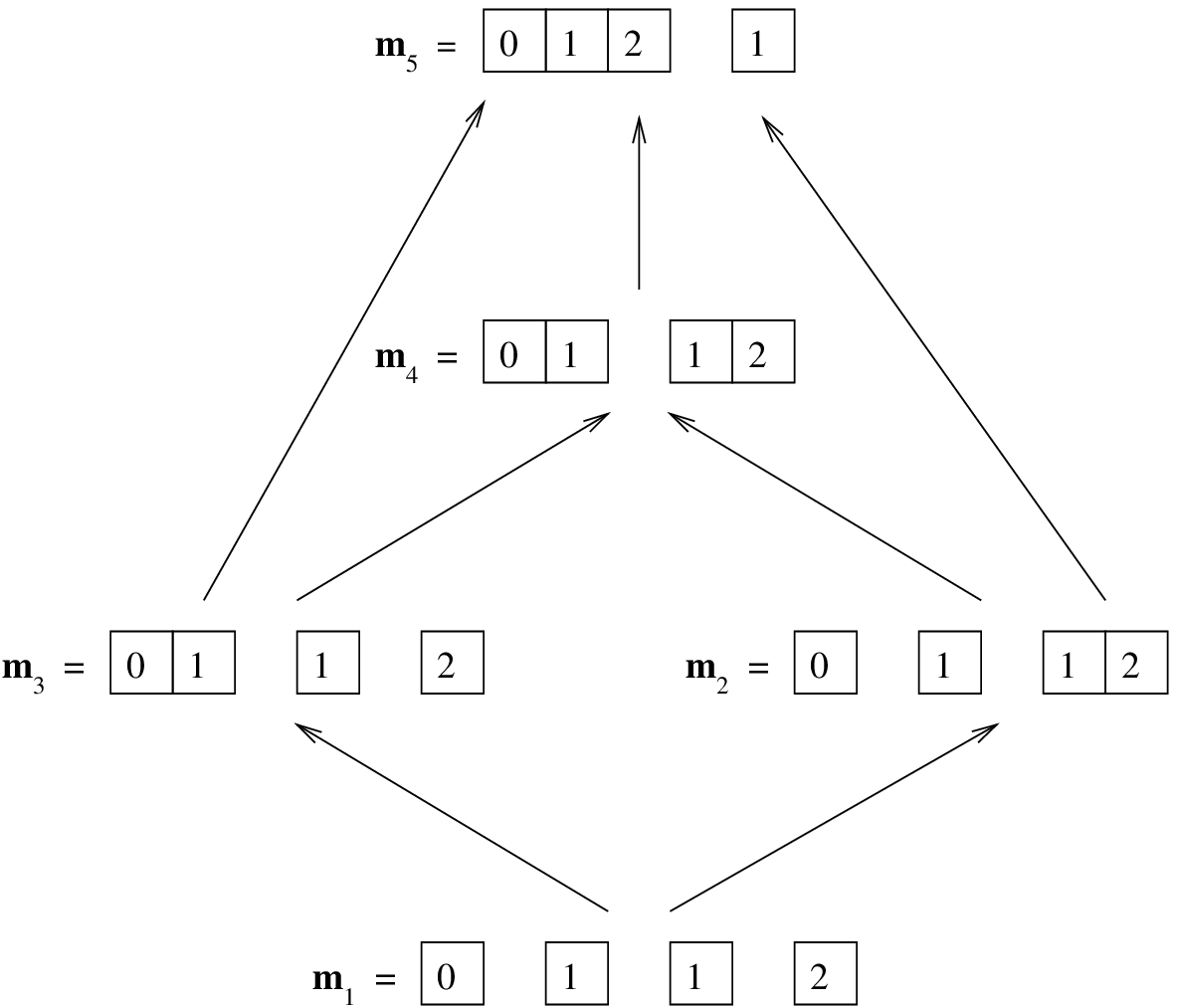}
\end{center}
\caption{\label{FIG0} The multi-segments larger than $[0]+2[1]+[2]$}
\end{figure}
\begin{example}\label{EX0}{\rm Let $\m_1 = [0] + 2[1] + [2]$.
This a minimal element for $\unlhd$. The multi-segments larger than
$\m_1$ for $\unlhd$ are 
\[
\m_2=[0]+[1]+[1,2],\ \m_3=[0,1]+[1]+[2],\ 
\m_4=[0,1]+[1,2],\ \m_5 =[1]+[0,2],
\]
and the poset structure is illustrated in Figure~\ref{FIG0}.
One has (\cite{Zel80}, 11.4):
\[
\left\{
\begin{array}{rcl}
[M_{\m_5}] &=& [L_{\m_5}] \\[1mm]
[M_{\m_4}] &=& [L_{\m_4}] + [L_{\m_5}]   \\[1mm]
[M_{\m_3}] &=& [L_{\m_3}] + [L_{\m_4}] +[L_{\m_5}]  \\[1mm]
[M_{\m_2}] &=& [L_{\m_2}] + [L_{\m_4}] +[L_{\m_5}] \\[1mm]
[M_{\m_1}] &=& [L_{\m_1}] + [L_{\m_2}] +[L_{\m_3}]+ 2\,[L_{\m_4}]+[L_{\m_5}]\,,
\end{array}
\right.
\]
which, by solving the equations, gives:
\[
\left\{
\begin{array}{rcl}
[L_{\m_5}] &=& [M_{\m_5}] \\[1mm]
[L_{\m_4}] &=& [M_{\m_4}] -[M_{\m_5}]   \\[1mm]
[L_{\m_3}] &=& [M_{\m_3}] -[M_{\m_4}]  \\[1mm]
[L_{\m_2}] &=& [M_{\m_2}] -[M_{\m_4}]\\ [1mm]
[L_{\m_1}] &=& [M_{\m_1}] -[M_{\m_2}] -[M_{\m_3}]+ [M_{\m_5}]\,.
\end{array}
\right.
\]
\finex
}
\end{example}

\subsection{} \label{SECT12BIS}
The composition multiplicities of the induced modules $M_\m$
are given by the so-called $p$-adic analogue of the Kazhdan-Lusztig 
conjecture formulated by Zelevinsky \cite{Zel81} and proved
by Ginzburg \cite{CG}, Theorem 8.6.23.
(See also \cite{Su} for another proof in the case of the degenerate
affine Hecke algebra.) 
Let us recall this result.

Fix a multi-segment $\m= \sum_{i\le j} m_{ij}[i,j]$,
and set $d_k = \sum_{i\le k \le j} m_{ij} \ (k\in\Z)$.
Clearly, only a finite number of $d_k$ are non-zero.
We call $d = (d_k)_{k\in\Z}$ the weight of $\m$
and write $\wt(\m) = d$. 

Let $V=\bigoplus_{k\in\Z} V_k$ be a $\Z$-graded vector space over $\C$
with $\dim V_k = d_k$. 
Let $E_V$ be the set of endomorphisms $x$ of $V$ of degree $-1$,
\ie such that $xV_k \subset V_{k-1}$.
(Equivalently, if $\Gamma$ denotes the quiver of type $A_\infty$
with orientation $k \longrightarrow k-1$ for all $k\in\Z$, then
$E_V$ is the set of representations of $\Gamma$ over $\C$ 
whose graded dimension is $d$.)

The group $G_V = \prod_{k\in\Z} GL(V_k)$ acts on $E_V$ by conjugation.
(The orbits of this action are precisely the isomorphism classes
of representations of $\Gamma$ with dimension~$d$.)
Note that all $x\in E_V$ are nilpotent and admit 
a Jordan basis $b$ consisting of homogeneous elements. 
Define the graded Jordan type of $x\in E_V$ as the multi-segment
$\n = \sum_{i\le j} n_{ij} [i,j]$, where $n_{ij}$ is the number
of Jordan cells of $b$ starting in $V_j$ and ending in $V_i$. 
It is easy to see that the $G_V$-orbit of $x$ consists of 
those elements $y$ having the same graded Jordan type as $x$.
Hence, the $G_V$-orbits $O_\n$ in $E_V$ are parametrized by the  
multi-segments $\n$ of weight $d$.

Zelevinsky has shown \cite{Zel81} that 
the closure of $O_\n$ decomposes as
$\overline{O}_\n = \coprod_{\m\unlhd\n} O_\m$, so that, 
$\m \unlhd \n$ if and only if $O_\m \subset \overline{O}_\n$.
Let ${\cal H}^i(\overline{O}_\n)_\m$ denote the stalk at a point $x\in O_\m$
of the $i$th intersection cohomology sheaf of the variety $\overline{O}_\n$.
\begin{theorem}[\cite{Zel81}, \cite{CG}]\label{THGINZ}
The composition multiplicity of the simple module $L_\n$ in the induced
module $M_\m$ is equal to
\[
K_{\m\n} = \sum_{i\ge 0} \dim {\cal H}^i(\overline{O}_\n)_\m \,.
\]
\end{theorem}
In \cite{Zel85}, Zelevinsky has further shown that the varietes 
$\overline{O}_\n$
are locally isomorphic to some Schubert varieties of type $A_{m-1}$,
where $m=\deg(\n)$.
Hence $K_{\m\n}$ is the value at $v=1$ of a certain Kazhdan-Lusztig polynomial
$K_{\m\n}(v)$ for the symmetric group $\SG_m$.

\subsection{} \label{SECT13}
Let $H_m$ denote the subalgebra of $\H_m$ generated by the elements
$T_i \ (1\le i \le m-1)$.
This is a Iwahori-Hecke algebra of type $A_{m-1}$.
The irreducible $H_m$-modules  are known to be
parametrized by the partitions $\alpha$ of $m$.
Let us denote them by $S_\alpha$.
In particular, $S_{(m)}$ and $S_{(1^m)}$ stand for the one-dimensional
$H_m$-modules associated to the characters $T_i \mapsto t$ and
$T_i \mapsto -1$, respectively.

For $z\in\C^*$ there exists a unique surjective homomorphism 
$\tau_z : \H_m \rightarrow H_m$ such that 
\[
\tau_z(T_i) = T_i\quad (1\le i \le m-1), \qquad\tau_z(y_1)=z.
\]
This is called the evaluation of $\H_m$ at $z$.
By pulling back the $H_m$-module $S_\alpha$ via $\tau_z$,
one gets an irreducible representation of $\H_m$ denoted
by $S_\alpha(z)$ and called an evaluation module. 
It belongs to $\CC_m^\Z$ when $z=t^a$ for some $a\in\Z$.
Consider the Young diagram $Y_\alpha(a)$ of $\alpha$ with each cell $c$
containing the integer $j-i+a$, where $i$ (\resp $j$) is the row index
(\resp column index) of $c$. 
In other words, the cells of $Y_\alpha(a)$ contain their content
shifted by $a$.
Each row of $Y_\alpha(a)$ can be regarded in a natural way as a segment,
and the collection of rows of $Y_\alpha(a)$ thus gives rise to a 
multi-segment 
\[
\m(\alpha,a) = \sum_{i=1}^r [a-i+1,a-i+\alpha_i] \,.
\]
Then $S_\alpha(t^a)$ is isomorphic
to $L_{\m(\alpha,a)}$ \cite{Ch}.

\begin{example}{\rm Among the irreducible modules of Example~\ref{EX0}, 
the only evaluation module is $L_{\m_4}$ which is isomorphic to
$S_{(2,2)}(t)$.
}
\finex
\end{example}

\subsection{}
Given a decomposition $m =m_1 + m_2$ of $m$ and some 
$\H_{m_i}$-modules
$M_i \ (i=1,2)$, one can form the induced module
\[
M_1 \odot M_2 := M_1 \otimes M_2 \uparrow_{\H_{m_1} \otimes \H_{m_2}}^{\H_m}.
\]
Here $\H_{m_1} \otimes \H_{m_2}$ is identified to the subalgebra
$\H_{(m_1,m_2)}$ of $\H_m$.
Suppose that the modules $M_i$ are simple, \ie $M_i = L_{\m_i}$
for some multi-segments $\m_i$. 
The problem we are interested in is whether one can formulate
in terms of $\m_1$ and $\m_2$ a necessary
and sufficient condition of irreducibility for the induced module
$M_1 \odot  M_2$. 

\subsection{}
Our strategy will be to reformulate this problem in the language
of canonical bases. 
Following Zelevinsky \cite{Zel80}, let us introduce
\[
R:=\bigoplus_{m\ge0} R_m\,,
\]
where we have put for convenience $R_0 = \C$.
By \ref{SECT12}, $\{[M_\m]\,|\, \m\in\M\}$ and $\{[L_\m]\,|\, \m\in\M\}$
are two bases of the vector space $R$.

The induction product endows $R$ with the structure of an associative
algebra:
\[
[M_1][M_2] := [M_1 \odot M_2].
\]
This algebra is commutative because $M_1 \odot M_2$ and $M_2 \odot M_1$
have the same composition factors with the same multiplicities
(\cite{Zel80}, Theorem 1.9).
Zelevinsky has shown that $R$ is in fact the polynomial ring 
in the variables $[L_\s] \ (\s\in \SS)$
(\cite{Zel80}, Corollary 7.5).
Using the restrictions from $\H_m$ to the subalgebras
$\H_{(k,m-k)}\cong \H_k\otimes\H_{m-k}$, one can define in a
standard way a comultiplication $c$ on $R$, and $R$ endowed with
these two operations becomes a graded bialgebra (\cite{Zel80}, Proposition
1.7), the grading being given by $\deg(R_m) = m$.
In particular $c$ is an algebra homomorphism determined by its
expression on the generators $[L_\s]$:  
\[
c[L_{[i,j]}] = 1\otimes[L_{[i,j]}] + 
\sum_{k=i}^{j-1} [L_{[i,k]}]\otimes [L_{[k+1,j]}] 
+[L_{[i,j]}]\otimes 1\,,
\]  
(\cite{Zel80}, Proposition 3.4).
Note that this formula shows that $c$ is not cocommutative.

On the other hand, let $N_\infty$ be the group
of upper triangular unipotent $(\Z\times \Z)$-matrices
with finitely many non-zero entries off the main diagonal.
We denote by $t_{ij}\ (i<j)$ the coordinate function
$(a_{kl})\in N_\infty \mapsto a_{ij}$.
Let $A$ be the ring of functions on $N_\infty$ which
are polynomials in the $t_{ij}$'s. 
The multiplication of $N_\infty$ induces a natural
comultiplication $\delta$ on $A$, given on the generators by
\[
\delta t_{ij} = 1\otimes t_{ij} +
\sum_{k=i+1}^{j-1} t_{ik}\otimes t_{kj}
+ t_{ij}\otimes 1\,,
\]
and this endows $A$ with the structure of a graded bialgebra,
where we put $\deg\, t_{ij} = j-i$.
The following simple but crucial observation immediately follows:
\begin{proposition} \label{ISO}
The graded bialgebras $R$ and $A$ are isomorphic via the map 
\[
\Psi[L_{[i,j]}] = t_{i,j+1}\,.
\]
In this isomorphism, the class 
$[M_\m]=\prod_{i\le j}[L_{[i,j]}]^{m_{ij}}$ of a  
generalized principal series representation is mapped to the monomial
$t_\m = \prod_{i\le j} t_{i,j+1}^{m_{ij}}$
in the coordinate functions.
\end{proposition}
A similar remark was made in \cite{Zel80}, 7.6.


\smallskip
In the next section, we shall see that $A$ has a natural quantum
deformation $A_v$, and that it allows us to define a canonical basis
of $A$ by specializing at $v=1$ a canonical basis of $A_v$.


\section{Quantum algebras and canonical bases} 
\label{SECT3}

We review the definitions of the quantum algebras $U^+_v$,
$A_v$ and of their canonical bases.
The main reference for this section are \cite{Lu90,BZ}. See also \cite{LTV}.

\subsection{}
Let $\nn_\infty$ denote the Lie algebra of 
strictly upper triangular $(\Z\times\Z)$-matrices with finitely many 
non-zero entries.
The standard basis of matrix units in $\nn_\infty$ will be denoted
by $\{e_{ij},\,i<j\in\Z\}$, and we write $e_i:=e_{i,i+1}$.
The enveloping algebra $U^+=U(\nn_\infty)$ is generated by the
$e_i,\ i\in\Z$, subject to Serre's relations:
\[
\begin{array}{ll}
e_ie_j = e_je_i, & |i-j|>1, \\[1mm]
e_i^2e_j -2e_ie_je_i+e_je_i^2 = 0, \qquad & |i-j|=1.
\end{array}
\]
Let $U^+_v$ be the associative $\Q(v)$-algebra
generated by elements $E_i \ (i\in\Z)$ subject to the relations
\[
\begin{array}{ll}
E_iE_j = E_jE_i, & |i-j|>1, \\[1mm]
E_i^2E_j -(v+v^{-1})E_iE_jE_i+E_jE_i^2 = 0, \qquad & |i-j|=1.
\end{array}
\]
This is the quantum enveloping algebra of $\nn_\infty$ in 
the sense of Drinfeld and Jimbo.
It can also be seen as the (twisted) Hall-Ringel algebra associated to 
an orientation $\Gamma$ of the Dynkin diagram of type $A_\infty$.
We shall always take for $\Gamma$ the standard orientation
$k \longrightarrow k-1$, $k\in\Z$ (\cf \S \ref{SECT12BIS}).

Let $\N^{(\Z)}$ be the semi-group of sequences $(d_j)_{j\in\Z}$
of non-negative integers with finitely many non-zero terms. 
Denote by $\alpha_i$ the sequence whose $i$th term is equal to 1
and all other terms are zero. 
We define a bilinear form on $\N^{(\Z)}$ by 
\[
(\alpha_i\,,\,\alpha_j) = \left\{
\begin{array}{cl}
2  & \mbox{if $i=j$,} \\[1mm]
-1 & \mbox{if $|i-j|=1$,} \\[1mm]
0  & \mbox{otherwise.}
\end{array}
\right.
\]
$\N^{(\Z)}$ identifies to the positive part of the root lattice
of the root system $A_\infty$, the $\alpha_i$ being the simple
roots, and the 
$\alpha_{ij}:=\alpha_i+\alpha_{i+1}+\cdots +\alpha_j$, $(i\le j)$,
the positive roots.
The algebra $U^+_v$ is $\N^{(\Z)}$-graded via the weight function
$\wt(E_i) = \alpha_i$.
The homogeneous components of $U^+_v$ are finite-dimensional,
and their linear bases are naturally labelled by multi-segments.
More precisely, the homogeneous component of weight $\alpha$ of $U^+_v$ has 
dimension equal to the value $p(\alpha)$ of the Kostant partition
function 
(this follows from the Poincar\'e-Birkhoff-Witt theorem for $U^+_v$). 
This value clearly coincides with 
the number of multi-segments $\m$ with $\wt(\m)=\alpha$. 

We shall also use the $\N$-grading of $U^+_v$ defined by 
$\deg(E_i)=1,\ i\in\Z$.

\subsection{} \label{SECT22}
Lusztig \cite{Lu90} has defined certain bases of $U^+_v$ associated
to orientations of the Dynkin diagram, called PBW-bases 
(they specialize when $v\mapsto 1$ to  bases of $U^+$ 
of the type provided by the Poincar\'e-Birkhoff-Witt theorem). 
We shall introduce the PBW-basis corresponding to the quiver $\Gamma$.
Let us describe the vector $E(\m)$ of this basis labelled by the 
multi-segment $\m$.
When $\m = [i,j]$ is reduced to a single segment, then $E(\m)$
is simply an iterated $v$-bracket. Namely 
\[
E([i])= E_i, \qquad
E([i,j]) = [E_j,[E_{j-1},[\cdots,[E_{i+1},E_i]_v\cdots ]_v]_v]_v,
\]
where $[x,y]_v := xy-v^{-(\alpha,\beta)}yx$ for $x,y\in U^+_v$ with
$\wt(x) = \alpha$ and $\wt(y) = \beta$.
Thus 
\[
E([1,2]) = E_2E_1 - v E_1E_2,\qquad
E([1,3]) = E_3E_2E_1-vE_3E_1E_2 -v E_2E_1E_3 + v^2 E_1E_2E_3.
\]
Clearly, $E([i,j])$ is a $v$-analogue of the root vector
$(-1)^{j-i}e_{i,j+1}$ of $\nn_\infty^+$, and we have
$\wt(E([i,j])=\alpha_{ij}$.
Next we introduce a total order on the set $\SS$ of segments by 
\[
[i,j] < [k,l] \quad \Longleftrightarrow \quad
\left\{
\begin{array}{l}
j<l \\[1mm]
\mbox{or}\\[1mm]
\mbox{$j=l$ and $i<k$.}
\end{array}
\right.
\]
(This coincides with the total order on positive roots associated with
the quiver $\Gamma$.)
Then the element of the PBW-basis indexed by $\m= \sum_{i\le j} m_{ij}[i,j]$ is
\[
E(\m) = \overrightarrow{\prod_{[i,j]\in\SS}}\, {1\over [m_{ij}]_v!} \, 
E([i,j])^{m_{ij}} \,,
\]
where for $a\in \Z$, we set $[a]_v=(v^a-v^{-a})/(v-v^{-1})$ and
$[a]_v!=[a]_v[a-1]_v\cdots [2]_v$. 
The arrow indicates that the product is taken in the order $<$ on $\SS$.
\subsection{}\label{SECT23}
To define the canonical basis of $U_v^+$ we consider the $\Z[v]$-lattice
\[
\L := \bigoplus_{\m\in\M} \Z[v] \, E(\m) \subset U_v^+\,,
\]
and the involution $x \mapsto \bar{x}$, defined as the unique ring automorphism
of $U_v^+$ such that
\[
\bar{v} = v^{-1}, \qquad \bar{E_i} = E_i.
\]
Lusztig has shown \cite{Lu90} that there exists a unique $\Q(v)$-basis 
$\{G(\m) \, | \, \m\in\M \}$ of $U^+_v$ such that
\[
\bar{G(\m)} = G(\m), \qquad G(\m) \equiv E(\m) \mod v\L .
\]
This is Lusztig's canonical basis (or Kashiwara's lower global crystal basis).
\begin{example}\label{EX2}{\rm Consider the homogeneous component of $U_v^+$ of
degree $\alpha_0 +2 \alpha_1 + \alpha_2$. It has dimension 5, and its 
weight vectors are labelled by the multi-segments (\cf Example~\ref{EX0}):
\[
\begin{array}{l}
\m_1 = [0]+2[1]+[2],\ \m_2=[0]+[1]+[1,2],\ \m_3=[0,1]+[1]+[2],\\ 
\m_4=[0,1]+[1,2],\ \m_5 =[1]+[0,2].
\end{array}
\]
The expansion of the $G(\m_i)$'s on the PBW-basis is
\[
\left\{
\begin{array}{rcl}
G(\m_1) &=& E(\m_1) \\[1mm]
G(\m_2)   &=& E(\m_2) +v^2\,E(\m_1)   \\[1mm]
G(\m_3)   &=& E(\m_3) + v^2\,E(\m_1)  \\[1mm]
G(\m_4)   &=& E(\m_4) + v\,E(\m_3)+ v\,E(\m_2)+(v+v^3)\,E(\m_1)\\[1mm]
G(\m_5)   &=& E(\m_5) + v\,E(\m_4) +v^2\,E(\m_3)+ v^2\,E(\m_2)+v^4\,E(\m_1)
\end{array}
\right.
\]
It can be calculated by using, for example, the formulas of \cite{LTV}
for the products $E_i\,E(\m)$.
\finex
}
\end{example}
Lusztig has also given a geometrical description of the canonical
basis, in terms of the varieties $\overline{O}_\m$ introduced in
\ref{SECT12BIS}.
\begin{theorem}[\cite{Lu90}]\label{THLU}
The expansion of $G(\n)$ on the basis $\{E(\m)\}$ is given by
\[
G(\n) = \sum_{\m\unlhd\n} K_{\m\n}(v) \, E(\m)\,,
\]
where, 
\[
K_{\m\n}(v) = v^{\dim O_\n -\dim O_\m}\ \sum_{i\ge 0} v^{-i}\,
\dim {\cal H}^i(\overline{O}_\n)_\m \,.
\] 
\end{theorem}
In particular $E(\m)$ occurs with a nonzero coefficient in the
expansion of $G(\n)$ if and only if $\m \unlhd \n$
for the partial order $\unlhd$ defined in \ref{SECT12}. 
Moreover, the coefficient of $E(\m)$ in  $G(\m)$ is 1.

The fact that precisely the same graded nilpotent orbits occur in
Theorem~\ref{THGINZ} and Theorem~\ref{THLU} is the deep geometrical fact
tying together the ${\mathfrak p}$-adic groups $GL_m$ and the
quantum group $U_v^+$.

\subsection{}\label{SECT24}
The canonical basis is `almost orthonormal' with respect to a scalar
product introduced by Kashiwara, which on the PBW-basis is given by
\[
(E(\m) \,,\, E(\n)) =
{(1-v^2)^{\deg(\m)}
\over
\prod_{i\le j}
\varphi_{m_{ij}}(v^2) }
\, \delta_{\m\,,\, \n} \,,
\]
where for $k\in\N$ we set $\varphi_k(z) = (1-z)(1-z^2)\cdots(1-z^k)$
(see \cite{LTV}, 4.3). By almost orthonormal one means that
\[
(G(\m)\,,\,G(\n)) \equiv \delta_{\m\,,\, \n} \mod v\AA,
\]
where $\AA$ is the subring of $\Q(v)$ consisting of functions regular at $v=0$. 
We will denote by $\{E^*(\m)\}$ and $\{G^*(\m)\}$ the adjoint bases
of $\{E(\m)\}$ and $\{G(\m)\}$ with respect to this scalar product.
Since $\{E(\m)\}$ is orthogonal, we see that $\{E^*(\m)\}$ is simply
a rescaling of $\{E(\m)\}$, namely
\[
E^*(\m) =  {\prod_{i\le j} \varphi_{m_{ij}}(v^2) \over 
(1-v^2)^{\deg(\m)}}\, E(\m)
= \overrightarrow{\prod_{[i,j]\in\SS}}\, v^{m_{ij}\choose 2} \,
E^*([i,j])^{m_{ij}}.
\]
It follows from Theorem~\ref{THLU} that 
\begin{equation}\label{EQLU}
E^*(\m) = \sum_{\m\unlhd\n} K_{\m\n}(v) \, G^*(\n)\,.
\end{equation}
So if we know the expansion of the canonical basis $\{G(\m)\}$ 
on the PBW-basis,
we can obtain the expansion of the dual canonical basis $\{G^*(\m)\}$ on the 
dual PBW-basis by solving this triangular system of linear
equations. (For a better algorithm, see below Section~\ref{SECT4}.)
In particular, we see that $E^*(\n)$ occurs in the expansion
of $G^*(\m)$ only if $\m\unlhd \n$, and the
coefficient of $E^*(\m)$ in $G^*(\m)$ is 1.
Hence for a single segment $\s\in\SS$ we have
$G^*(\s) = E^*(\s)$.
\begin{example}\label{EX3}{\rm Retaining the notation of Example~\ref{EX2}, 
we have
\[
\left\{
\begin{array}{rcl}
E^*(\m_5) &=& G^*(\m_5) \\[1mm]
E^*(\m_4) &=& G^*(\m_4) +v\,G^*(\m_5)   \\[1mm]
E^*(\m_3) &=& G^*(\m_3) + v\,G^*(\m_4) +v^2\,G^*(\m_5)  \\[1mm]
E^*(\m_2) &=& G^*(\m_2) + v\,G^*(\m_4) +v^2\,G^*(\m_5) \\[1mm]
E^*(\m_1) &=& G^*(\m_1) + v^2\,G^*(\m_2) +v^2\,G^*(\m_3)+ (v+v^3)\,G^*(\m_4)
+v^4\,G^*(\m_5)\,,
\end{array}
\right.
\]
and
\[
\left\{
\begin{array}{rcl}
G^*(\m_5) &=& E^*(\m_5) \\[1mm]
G^*(\m_4) &=& E^*(\m_4) -v\,E^*(\m_5)   \\[1mm]
G^*(\m_3) &=& E^*(\m_3) -v\,E^*(\m_4)  \\[1mm]
G^*(\m_2) &=& E^*(\m_2) - v\,E^*(\m_4)\\ [1mm]
G^*(\m_1) &=& E^*(\m_1) - v^2\,E^*(\m_2) -v^2\,E^*(\m_3)+ 
(v^3-v)\,E^*(\m_4)+v^2\,E^*(\m_5)\,.
\end{array}
\right.
\]
\finex
}
\end{example}

\subsection{}\label{SECT25}
Let $A_v$ be the quantum analogue of the algebra $A$ of polynomials
in the coordinate functions $t_{ij}$  of the group $N_\infty$.
We denote by $T_{ij}\ (i<j)$ the $v$-analogue of 
$t_{ij}$ 
(there should be no risk of confusing the $T_{ij}$'s with 
the generators $T_i$ of $\H_m$).
It is often convenient to write $T_{ij}=0$ for $i>j$, $T_{ii}=1$, and
to index the non-trivial $T_{ij}$'s by segments, that is,
$T_{ij}=T_{[i,j-1]}$, for $i< j$.

The commutation relations satisfied by the $T_\s$, $(\s\in \SS)$ 
are the following
(\cite{BZ}, Proposition 3.11).
Let $\s$ and $\s'$ be two segments such that $\s'>\s$
for the total ordering defined in \ref{SECT22}. Then
\begin{equation}\label{REDRESS}
T_{\s'}T_\s =
\left\{
\begin{array}{ll}
v^{-(\wt(\s'),\wt(\s))}\,\left(T_\s T_{\s'} + 
(v^{-1}-v)T_{\t'}T_{\t}\right)
&\mbox{if $\s$ and $\s'$ are linked,} \\[1mm]
v^{-(\wt(\s'),\wt(\s))}\,T_\s T_{\s'} &\mbox{otherwise.}
\end{array}
\right.
\end{equation}
As in \ref{SECT12},
in the case where $\s$, $\s'$ are linked, we have put 
$\t=\s\cup \s'$ and $\t'=\s\cap \s'$. When $\t' = \emptyset$, 
we understand $T_{\t'}=1$.

In fact, as proved in \cite{BZ}, the algebras $U_v^+$ and $A_v$
are isomorphic, the isomorphism $\Phi$ being given by 
$\Phi(E_i) = T_{[i]}=T_{i,i+1}$.
More generally, we have $\Phi(E^*([i,j])) = T_{[i,j]}$ 
(which motivates the notation introduced above), so that
\[
\Phi(E^*(\m)) = \overrightarrow{\prod_{[i,j]\in\SS}}\, v^{m_{ij}\choose 2} \,
T_{[i,j]}^{m_{ij}}. 
\]
For $\m\in\M$, we shall write 
\[
T_\m =  \overrightarrow{\prod_{[i,j]\in\SS}} T_{[i,j]}^{m_{ij}}.
\]
Relations~(\ref{REDRESS}) can be regarded as a set of straightening 
rules for computing the expansion of an arbitrary polynomial
in the $T_{ij}$'s on the basis $\{T_\m\}$.

\subsection{}\label{SECT26}
Although the algebras $U_v^+$ and $A_v$ are isomorphic, their natural
specializations at $v=1$ are not. 
Indeed, to specialize $U_v^+$ one first considers the integral form 
\[
U_{v,\Z}^+=\bigoplus_{\m\in\M} \Z[v,v^{-1}]\,E(\m),
\] 
and then one sets
$U_1^+ = \C\otimes_{\Z[v,v^{-1}]}U_{v,\Z}^+$, where $\C$ is regarded
as a $\Z[v,v^{-1}]$-module via $v\mapsto 1$.
Clearly, $U_1^+\simeq U^+$ via $1\otimes E_i \mapsto e_i$.
Whereas for $A_v$ one takes 
\[
A_{v,\Z}=\bigoplus_{\m\in\M} \Z[v,v^{-1}]\,T_\m,
\] 
and then $A_1=\C\otimes_{\Z[v,v^{-1}]}A_{v,\Z}$, which
is isomorphic to $A$ by $1\otimes T_{ij} \mapsto t_{ij}$.
Thus, $U_v^+$ specializes to the non-commutative enveloping 
algebra $U^+$,
while $A_v$ specializes to the commutative algebra of polynomial 
functions $A$. 
\begin{example}{\rm In $U_v^+$ we have 
\begin{equation}\label{STAR}
E([1,2]) = E_2E_1 - v E_1E_2
\end{equation}
which gives in $U^+$ the familiar relation 
$
-e_{13}=e_{23}e_{12}-e_{12}e_{23}.
$
Relation (\ref{STAR}) is transformed under $\Phi$ into 
\[
(1-v^2)T_{[1,2]} = T_{[2]}T_{[1]}-vT_{[1]}T_{[2]}
\]
which specializes in $A$ to $t_{23}t_{12}=t_{12}t_{23}$.
}
\finex
\end{example}
In the sequel the algebras $U_v^+$ and $A_v$ will in general
be identified via $\Phi$,
and we shall distinguish between them only when the specialization 
$v\mapsto 1$ is considered.
In particular, the basis $\{G^*(\m)\}$ of $U^+_v$ gets identified 
under $\Phi$ to the dual canonical
basis (or string basis) of $A_v$ studied by Berenstein and
Zelevinsky in \cite{BZ,BZ2}.

\subsection{}

The bases $\{E(\m)\}$ and $\{G(\m)\}$ 
(\resp $\{E^*(\m)\}$ and $\{G^*(\m)\}$) give rise at $v=1$ to
well-defined bases of $U^+$ (\resp $A$).
We  denote them by  $\{e(\m)\}$, $\{g(\m)\}$,
$\{e^*(\m)\}$ and $\{g^*(\m)\}$, respectively.
We put $\B_v := \{G^*(\m)\}$ and $\B := \{g^*(\m)\}$. 
It follows from Equation~(\ref{EQLU}) that
\begin{equation}\label{EQLUSP}
e^*(\m) = \sum_{\m\unlhd\n} K_{\m\n} \, g^*(\n)\,.
\end{equation}
Hence, comparing with Theorem~\ref{THGINZ}, we obtain 
\begin{theorem}\label{THAR}
Under the isomorphism $\Psi : R \longrightarrow A$ of
Proposition~{\rm\ref{ISO}}, the class of the 
simple module $L_\m$ is mapped to the element $g^*_\m$ of $\B$.
\end{theorem}
We note that Ariki has proved a much more general result \cite{Ar}. 
It describes in a similar way the simple modules of an infinite family of 
finite-dimensional quotients of $\H_m$ (the so-called 
cyclotomic Hecke algebras), and it also handles the case when
$t$ is an $\ell$th root of unity. For an application to the $\ell$-modular
representation theory of $GL_m$, see \cite{LTV}.

Proposition~\ref{ISO} and Theorem~\ref{THAR} imply immediately 
the following criterion of 
irreducibility for an induction product of simple modules.
\begin{theorem}\label{TH10}
Let $\m,\,\n$ be two multi-segments of degree $m$ and $n$ respectively,
and $L_\m,\,L_\n$ be the corresponding simple modules over $\H_m$ and
$\H_n$. 
The induction product $L_\m\odot L_\n$ is a simple $\H_{m+n}$-module
if and only if the product $g^*(\m)\,g^*(\n)$ belongs to $\B$.
\end{theorem}
\begin{example}{\rm Take $\m=[1]+[2,3]$ and $\n=[2]+[3,4]$.
Then, using for example the algorithm described in the next
section, one calculates
\[
G^*(\m)\,G^*(\n) = 
v^{-1}\,G^*(\p_1)
+G^*(\p_2)
+G^*(\p_3)
+v\,G^*(\p_4)
+G^*(\p_5),
\]
where $\p_1=[1]+[2]+[2,3]+[3,4]$,
$\p_2=[1]+[2]+[3]+[2,4]$,
$\p_3=[1,2]+[2,3]+[3,4]$,
$\p_4=[1,2]+[3]+[2,4]$,
$\p_5=[1,3]+[2,4]$.
Therefore
\[
g^*(\m)\,g^*(\n) = 
g^*(\p_1)
+g^*(\p_2)
+g^*(\p_3)
+g^*(\p_4)
+g^*(\p_5),
\]
hence $L_\m\odot L_\n$  is not irreducible 
(which agrees with 
Corollary~\ref{COR27} below) and its composition
factors are $L_{\p_1}, \ldots , L_{\p_5}$, each of them appearing with
multiplicity 1.
}
\finex
\end{example}

\begin{proposition}\label{PROP15}
The product $g^*(\m)\,g^*(\n)$ belongs to $\B$ if and only if there
exists $k \in \Z$ such that $v^k\,G^*(\m)\,G^*(\n)$ belongs to $\B_v$.
\end{proposition}
\proof
Given three multi-segments $\m,\,\n,\,\p$, define the coefficient
$\alpha_{\m\n}^\p(v)$ by
\begin{equation}\label{DEFALPHA}
G^*(\m)\,G^*(\n) = \sum_\p \alpha_{\m\n}^\p(v) G^*(\p)\,.
\end{equation}
Using his geometrical description of the canonical basis, 
Lusztig (\cite{Lu91}, 11.5) has proved the positivity of the
$\alpha_{\m\n}^\p(v)$,
namely
\begin{equation}\label{POSLU}
\alpha_{\m\n}^\p(v) \in \N[v,v^{-1}]\,.
\end{equation}
Now 
\begin{equation}
g^*(\m)\,g^*(\n) = \sum_\p \alpha_{\m\n}^\p(1) g^*(\p)\,,
\end{equation}
so if $g^*(\m)\,g^*(\n) = g^*(\q)$, then 
$\alpha_{\m\n}^\p(1) =\delta_{\p\q}$ and all $\alpha_{\m\n}^\p(v)$
are zero except $\alpha_{\m\n}^\q(v)$ which can only be a power of
$v$.
The converse is obvious.
\cqfd


\section{Calculation of the dual canonical basis} 
\label{SECT4}
\subsection{}
The basis $\{G^*(\m)\}$ can be characterized by two conditions similar 
to those defining $\{G(\m)\}$.
Let 
\[
\L^*:=\bigoplus_{\m\in\M} \Z[v] \, E^*(\m) \,,
\]
and let $\tau$ denote the anti-automorphism of $U_v^+$ such that
$\tau(E_i)=E_i$.
\begin{proposition}\label{PROP4}
Let $\m\in\M$ and write $|\m|^2:=(\wt(\m),\wt(\m))$. Then 
$G^*(\m)$ is the unique homogeneous element of degree $\wt(\m)$
of $U_v^+$ satisfying 
\[
\bar{G^*(\m)} = v^{\deg(\m)-{1\over 2}|\m|^2}\,\tau(G^*(\m)),\qquad
G^*(\m) \equiv E^*(\m) \mod v\L^* .
\]
\end{proposition}
\proof
That $G^*(\m)$ satisfies the second condition is obvious from the
definitions.
That it also satisfies the first condition was shown by Reineke
\cite{Rei}, Lemma 4.3. 
(Our power of $v$ is different because we 
use a different normalization of the scalar product.)
Let us prove unicity.
Suppose that $x\in U_v^+$ has weight $\wt(\m)$ and satisfies
$x\in E^*(\m)+v\L^*$.
Then
\[
x=G^*(\m)+
\sum_{\wt(\n)=\wt(\m),} 
\alpha_\n(v)\,G^*(\n)
\]
for some $\alpha_\n(v) \in v\Z[v]$.
If moreover $\bar{x} =  v^{\deg(\m)-{1\over 2}|\m|^2}\,\tau(x)$
then $\alpha_\n(v^{-1})=\alpha_\n(v),$ which forces
$\alpha_\n(v)=0$.
\cqfd

The following relation satisfied by the coefficients
$\alpha_{\m\n}^\p(v)$
defined in (\ref{DEFALPHA}) will be important in the
sequel.
\begin{corollary}[\cite{Rei}]\label{EQREI}
For $\m,\n,\p \in \M$,
\[
\alpha_{\n \m}^\p(v)=v^{-(\wt(\m),\wt(\n))}\,
\alpha_{\m \n}^\p(v^{-1})\,.
\]
\end{corollary}
\proof
On the one hand,
\begin{eqnarray*}
\bar{G^*(\m)\,G^*(\n)} &=&
v^{\deg(\m)+\deg(\n)-1/2(|\m|^2+|\n|^2)}
\,\tau(G^*(\m))\,\tau(G^*(\n))\\[1mm]
&=&v^{\deg(\m)+\deg(\n)-1/2(|\m|^2+|\n|^2)}
\,\tau(G^*(\n)\,G^*(\m)) \\[1mm]
&=&v^{\deg(\m)+\deg(\n)-1/2(|\m|^2+|\n|^2)}
\,\sum_\p \alpha_{\n\m}^\p(v) \,\tau(G^*(\p))\,. 
\end{eqnarray*}
On the other hand,
\begin{eqnarray*}
\bar{G^*(\m)\,G^*(\n)} &=&
\sum_\p \alpha_{\m\n}^\p(v^{-1})\, \bar{G^*(\p)} \\[1mm]
&=& 
\sum_\p \alpha_{\m\n}^\p(v^{-1})\, 
v^{\deg(\p)-1/2|\p|^2} \, \tau(G^*(\p))\,,
\end{eqnarray*}
and the result follows by comparison, since 
$\wt(\p) = \wt(\m) +\wt(\n)$ for any non-zero coefficient
$\alpha_{\m\n}^\p(v)$.
\cqfd

\subsection{}\label{SECT211}
Let $\Rg$ be the sub-ring of $\Z[v,v^{-1}]$ consisting of the 
Laurent polynomials invariant under $v\mapsto v^{-1}$.
We set $\V := \bigoplus_{\m\in\M}\Rg\, G^*(\m)$.
Then, by Proposition~\ref{PROP4}, $\V$ may be described as the 
$\Rg$-lattice whose homogeneous elements $x$ satisfy
\[
\bar x = v^{\deg(x)-{1\over 2}|\wt(x)|^2}\,\tau(x) \,,
\]
and the vector $G^*(\m)$ is characterized by
\[
G^*(\m) \in \V, \qquad G^*(\m) \equiv E^*(\m) \mod v\L^* .
\]

Let us describe an algorithm similar to that of \cite{LLT96}
for computing the transition matrix from $\{E^*(\m)\}$ to
$\{G^*(\m)\}$.
Given two multi-segments 
$\m=\sum_{\s\in\SS} m_\s\,\s$ and  $\n=\sum_{\s\in\SS} n_\s\,\s$,
define
\begin{eqnarray*}
b(\m,\n) &:=& \sum_{\s'>\s}
m_{\s'}n_{\s}(\wt(\s),\wt(\s')) + \sum_\s m_\s\,n_\s\\[1mm]
&=& \sum_{\s'>\s}
m_{\s'}n_{\s}(\wt(\s),\wt(\s')) + 
{1\over 2}\sum_\s m_\s\,n_\s(\wt(\s),\wt(\s)).
\end{eqnarray*}
Clearly we have
\begin{equation}\label{bmn}
b(\m,\n) + b(\n,\m) = (\wt(\m) \,,\,\wt(\n))\,.
\end{equation}
It follows from the straightening relations~(\ref{REDRESS}) satisfied by the 
elements $T_\s, \ \s \in \SS$ that
\[
E^*(\m)\,E^*(\n) = v^{-b(\m,\n)}\,E^*(\m+\n) + \mbox{higher terms}\,,
\]
where by `higher terms' we mean a linear combination of $E^*(\p)$
with $\m+\n \lhd \p$.
We know that $E^*(\m')$ occurs in the expansion of $G^*(\m)$ only if
$\m\unlhd\m'$.
It follows that $\alpha_{\m\n}^\p(v) \not = 0$ only if
$\m + \n \unlhd \p$, and for $\p = \m+\n$,
\begin{equation}\label{HIT}
\alpha_{\m\n}^{\m+\n} = v^{-b(\m,\n)}.
\end{equation}
Reineke \cite{Rei} has given an interesting representation-theoretical 
interpretation
of $b(\m,\n)$, namely, if $M$ and $N$ denote respectively the
representations of the quiver $\Gamma$ of isomorphism type $\m$ and
$\n$, then
\[
b(\m,\n) = \dim\hom(N,M) - \dim\ext(M,N)\,.
\] 
Now, introduce the following element of $U^+_v$:
\[
U(\m,\n):={1\over v-v^{-1}}\left(
v^{b(\m,\n)+1}\,G^*(\m)\,G^*(\n) - v^{b(\n,\m)-1}\,
G^*(\n)\,G^*(\m)
\right).
\]
\begin{proposition}\label{PROP13}
Writing $U(\m,\n)=\sum_\p \beta_{\m,\n}^\p(v)\,G^*(\p)$,
we have

\noindent
\begin{tabular}{ll}
{\rm (i)} &$\beta_{\m\n}^\p(v) \in \Z[v,v^{-1}]$ and 
$\beta_{\m\n}^{\m+\n}(v) = 1$,\\[1mm]
{\rm (ii)}& $\beta_{\m\n}^\p(v) \not = 0$ only if $\m+\n \unlhd \p$,\\[1mm]
{\rm (iii)} & $\beta_{\m\n}^\p(v^{-1}) = \beta_{\m\n}^\p(v)$.
\end{tabular}
\end{proposition}
\proof
Everything follows immediately from Corollary~\ref{EQREI} and
relations (\ref{bmn}) and (\ref{HIT}). 
\cqfd

Note that if 
\[
G^*(\m)\,G^*(\n) = v^{c(\m,\n)}\,G^*(\n)\,G^*(\m),
\]
then it follows from (\ref{HIT}) that 
$ 
c(\m,\n) = b(\m,\n)-b(\m,\n).
$
Hence in such a case $U(\m,\n)$ reduces to 
\[
U(\m,\n) = v^{b(\m,\n)}\,G^*(\m)\,G^*(\n)\,.
\]

We can now easily compute the expansion of the dual canonical basis 
$\{G^*(\n)\}$ on the dual PBW-basis
by induction on $\deg(\n)$. 
Suppose that $\{G^*(\n)\}$ has already been calculated for
$\deg(\n) \le n$.
Let $\m$ be a multi-segment of degree $n+1$.
If $\m$ consists of a single segment, we know that $G^*(\m) = E^*(\m)$.
If not, write $\m = \m_1 + \m_2$ for some non-empty
multi-segments $\m_1$ and $\m_2$. 
(For example, one can take for
$\m_2$ the largest segment occuring in $\m$ with a non-zero
multiplicity).
Then the expansions of $G^*(\m_1)$ and $G^*(\m_2)$ on the dual PBW-basis
are known by induction, and one can calculate the expansion 
of $V(\m) := U(\m_1,\m_2)$ on $\{E^*(\n)\}$ by means of the
straightening rules of \S\,\ref{SECT25}.
If $\m$ is a single segment, we set $V(\m) = E^*(\m)$.
 
Then, by Proposition~\ref{PROP13} (i) (ii), the transition matrices 
from 
$\{G^*(\m)\,|\, \deg(\m) = n+1\}$ to $\{V(\m)\,|\, \deg(\m) =n+1\}$ 
and from 
$\{E^*(\m)\,|\, \deg(\m) = n+1\}$ to $\{V(\m)\,|\, \deg(\m) = n+1\}$
will be both unitriangular if their
rows and columns are arranged in any total order $\le$ extending
the partial order $\unlhd$, and their entries will belong 
to $\Z[v,v^{-1}]$.
Moreover,  by Proposition~\ref{PROP13} (iii), $V(\m) \in \V$.
This implies that $\{V(\m)\,|\, \deg(\m) = n+1\}$ is 
an $\Rg$-basis of the degree $n+1$ homogeneous component of $\V$. 

Let $\{ \m_1, \ldots , \m_r \}$ be the list of all multi-segments
of degree $n+1$ and of a given weight $\lambda \in \N^{(\Z)}$,
arranged in increasing order with respect to $\le$. 
By Proposition~\ref{PROP13}~(i) and (ii), $V(\m_r)  = G^*(\m_r)$.
Assume by descending induction on $k$ that the expansions of 
$G^*(\m_{k+1}), \ldots , G^*(\m_r)$ on $\{E^*(\n)\}$ are known.
Then 
\[
G^*(\m_k) = V(\m_k) - \sum_{i=k+1}^r \gamma_i(v)\, G^*(\m_i)
\]
where the coefficients $\gamma_i(v)$ are completely determined
by the conditions
\[
\gamma_i(v^{-1}) = \gamma_i(v)\,,\qquad
G^*(\m_k) \equiv E^*(\m_k) \mod v\L^*\,.
\]
Hence we can obtain the basis $\{G^*(\m)\}$ from the previously
known basis $\{V(\m)\}$ by a simple algorithm.

\begin{example}{\rm Let us calculate with this algorithm the
vectors $G^*(\m_i)$ of Example~\ref{EX3}.
For $i=1,\ldots ,5$ we write $\m_i = \n_i + \s_i$ where $\s_i$
is the largest segment of $\m_i$.
Then, by induction we know that
\[
\left\{
\begin{array}{rcl}
G^*(\n_5) &=& G^*([1]) = E^*([1]) \\[1mm]
G^*(\n_4) &=& G^*([0,1]) = E^*([0,1])  \\[1mm]
G^*(\n_3) &=& G^*([0,1]+[1]) = E^*([0,1]+[1])  \\[1mm]
G^*(\n_2) &=& G^*([0]+[1]) = E^*([0]+[1]) - v\,E^*([0,1])\\ [1mm]
G^*(\n_1) &=& G^*([0]+2[1]) = E^*([0]+2[1]) - v^2\,E^*([0,1]+[1]) \,.
\end{array}
\right.
\]
On the other hand, since $\s_i$ is a segment, $G^*(\s_i) = E^*(\s_i)$.
The vectors $V(\m_i)=U(\n_i,\s_i)$ are then calculated and found to
be equal to
\[
\left\{
\begin{array}{rcl}
V(\m_5) &=& E^*(\m_5) \\[1mm]
V(\m_4) &=& E^*(\m_4) +v^{-1}\,E^*(\m_5)   \\[1mm]
V(\m_3) &=& E^*(\m_3) +v^{-1}\,E^*(\m_4)+v^{-2}\,E^*(\m_5)  \\[1mm]
V(\m_2) &=& E^*(\m_2) - v\,E^*(\m_4)\\ [1mm]
V(\m_1) &=& E^*(\m_1) +(1+ v^{-2})\,E^*(\m_2) -v^2\,E^*(\m_3) 
-v\,E^*(\m_4)-\,E^*(\m_5)\,.
\end{array}
\right.
\]
We see that $V(\m_2)$ and $V(\m_5)$ belong to $\B_v$, but the other
$V(\m_i)$'s have to be `corrected' since they do not specialize to
$E^*(\m_i)$ when $v \mapsto 0$. This gives
\[
\left\{
\begin{array}{rcl}
G^*(\m_5) &=& V(\m_5) = E^*(\m_5) \\[1mm]
G^*(\m_4) &=& V(\m_4) - (v+v^{-1})\, G^*(\m_5) 
              = E^*(\m_4) -v\,E^*(\m_5)   \\[1mm]
G^*(\m_3) &=& V(\m_3)-(v+v^{-1})\,G^*(\m_4)
              -(v^2+1+v^{-2})\,G^*(\m_5) \\[1mm]
          &=& E^*(\m_3) -v\,E^*(\m_4)  \\[1mm]
G^*(\m_2) &=&  V(\m_2) = E^*(\m_2) - v\,E^*(\m_4)\\ [1mm]
G^*(\m_1) &=& V(\m_1) - (v^2+1+v^{-2})\,G^*(\m_2)
                -(v+v^{-1})\,G^*(\m_4)  \\[1mm]
          &=& E^*(\m_1) - v^2\,E^*(\m_2) -v^2\,E^*(\m_3)+ 
(v^3-v)\,E^*(\m_4)+v^2\,E^*(\m_5)\,.
\end{array}
\right.
\] 
Note that there are many possible choices for $\{V(\m)\}$,
some of them giving better approximations of $\{G^*(\m)\}$ and
therefore more efficient algorithms. 
For example 
\[
V'(\m_3)=U([0,1],[1]+[2])= E^*(\m_3) -v\,E^*(\m_4) = G^*(\m_3)
\]
would be a better choice than $V(\m_3)$.
We shall not discuss this question here.
}
\finex
\end{example}


\section{The Berenstein-Zelevinsky conjecture} \label{SECT5} 

\subsection{}
The basis $\B_v$ seems to enjoy remarkable multiplicative properties.
Let us say that $x$ and $y$ {\em quasi-commute} if $xy=v^n\,yx$ for some $n\in\Z$.
The following deep conjecture was formulated by Berenstein and Zelevinsky,
together with some supporting evidence. 
\begin{conjecture}[\cite{BZ}]\label{CONJ2}
Let $\m,\,\n \in \M$. The product $G^*(\m)\,G^*(\n)$ belongs to $\B_v$
up to some power of $v$ if and only if $G^*(\m)$ and $G^*(\n)$ quasi-commute. 
\end{conjecture}
Clearly, Conjecture~\ref{CONJ2}  
would imply that more generally the product 
$G^*(\m_1)\cdots G^*(\m_r)$ belongs to $\B_v$
up to a power of $v$ if and only if $G^*(\m_k)$ and $G^*(\m_l)$
quasi-commute for all $1\le k<l \le r$.
\begin{example}\label{GPS}
{\rm For a segment $\s$ we have $G^*(\s) = E^*(\s)$. Hence, by 
Equation~(\ref{EQLU}), we deduce that if 
$\m = \s_1 + \cdots + \s_r$, the product 
$G^*(\s_1)\cdots G^*(\s_r)$ belongs to $\B_v$ up to some
power of $v$ if and only if $\m$ is maximal for the partial order
$\unlhd$,
that is, if and only if the $\s_i$ are pairwise not linked.
On the other hand, by Equation~(\ref{REDRESS}), this is equivalent
to the fact that the $G^*(\s_i)$ pairwise $v$-commute. 
}
\finex
\end{example}
One implication of Conjecture~\ref{CONJ2} readily follows 
from Corollary~\ref{EQREI} \cite{Rei}.
Indeed if 
$G^*(\m)\,G^*(\n) = v^k\,G^*(\p)$ for some $k\in\Z$ and
some $\p\in\M$, then 
$G^*(\n)\,G^*(\m) = v^{-k-(\wt(\m),\wt(\n))}\,\,G^*(\p)$,
which means that $G^*(\m)$ and $G^*(\n)$ quasi-commute.
Another way of seeing this is by specializing $v$ to $1$.
For if $v^k\,G^*(\p) = G^*(\m)\,G^*(\n)$, then 
$g^*(\p) = g^*(\m)\,g^*(\n) = g^*(\n)\,g^*(\m)$, and
because of the positivity~(\ref{POSLU}) we have
$G^*(\n)\,G^*(\m) = v^l\,G^*(\p)$ for some $l\in\Z$.

Recall also from \S \ref{SECT211} that if 
$G^*(\m)$ and $G^*(\n)$ quasi-commute, then 
\[
v^{b(\m,\n)}\,G^*(\m)\,G^*(\n)=U(\m,\n),
\]
and thus $v^{b(\m,\n)}\,G^*(\m)\,G^*(\n) \in \V$.

In view of these remarks, Conjecture~\ref{CONJ2}
may be reduced to the following
\begin{conjecture}\label{CONJ2'}
Let $\m$, $\n \in \M$. If $G^*(\m)$ and $G^*(\n)$ quasi-commute, then
\[
v^{b(\m,\n)}\,G^*(\m)\,G^*(\n) \equiv E^*(\m+\n) \mod v\L^*.
\]
\end{conjecture}

\begin{example}{\rm (Compare with Lemma~II.8 of \cite{MW}.) 
Let $b<b''<e'<e<e''$ be integers and consider
the segments 
\[
\s=[b,e],\qquad \s'=[b,e'],\qquad \s''=[b'',e'']\,.
\]
We have $\s' < \s < \s''$, and
\[
G^*(\s')=E^*(\s'),\quad 
G^*(\s+\s'')=E^*(\s+\s'')-v\,E^*(\t+\t'),
\]
where $\t = \s\cup \s''=[b,e'']$ and $\t'=\s\cap \s''=[b'',e]$.
From this, we readily obtain that
\[
G^*(\s')\,G^*(\s+\s'') = E^*(\s+\s'+\s'') -v\,E^*(\s'+\t + \t')
= v\,G^*(\s+\s'')\,G^*(\s').
\]
Thus, in this case Conjecture~\ref{CONJ2'} is verified and 
$$G^*(\s')\,G^*(\s+\s'') = G(\s'+\s+\s'')\,.$$
Hence 
$L_{\s'}\odot L_{\s+\s''}$ is irreducible and isomorphic to
$L_{\s'+\s+\s''}$.  
}
\finex
\end{example}

\subsection{} \label{SECT29}
For $k\in\Z$ and $r\in\N^*$, let $Q_{k,r}$ denote the sub-semigroup
of $\N^{(\Z)}$ generated by the simple roots 
$\alpha_k,\,\alpha_{k+1},\,\ldots,\alpha_{k+r-1}$.
If $\wt(\m)$ and $\wt(\n)$ belong to $Q_{k,r}$ for some $k\in\Z$ and $r\le 3$,
then Conjecture~\ref{CONJ2'} is true \cite{BZ}.
Moreover in this case the set
$\{G^*(\m)\ | \ \wt(\m) \in Q_{k,r}\}$
consists of quasi-commuting
products of a small number of certain special elements called 
quantum minors (4 primitive minors for $r=2$ and 12 for $r=3$). 
  
Given two subsets $I=\{i_1<\cdots <i_k\}$, 
$J=\{j_1<\cdots < j_k\}$ of $\Z$, 
the {\em quantum minor} $\Delta(I,J)$ is defined by
\[
\Delta(I,J):=\sum_{\sigma\in\SG_k}
(-v)^{\ell(\sigma)}\,T_{i_1j_{\sigma(1)}}\cdots T_{i_kj_{\sigma(k)}}.
\]
Note that with our convention for $T_{ij}$, $\Delta(I,J) \not = 0$ only if
$i_r\le j_r\ (r=1,\ldots ,k)$.
Moreover, if $i_r=j_r$ then $\Delta(I,J)$ factors into
\begin{equation}\label{EQQM}
\Delta(I,J) = \Delta(I',J')\,\Delta(I'',J'')
\end{equation}
where $I'=\{i_1,\ldots ,i_{r-1}\}$, $I''=\{i_{r+1},\ldots ,i_k\}$,
$J'=\{j_1,\ldots ,j_{r-1}\}$ and $J''=\{j_{r+1},\ldots ,j_k\}$.
One can thus assume that $i_r \le j_r\ (r=1,\ldots ,k)$ and
attach to $\Delta(I,J)$ the multi-segment 
\begin{equation}
\m(I,J)= \sum_{r=1}^k \,[i_r,j_r-1]
\end{equation}
in which $[i_r,j_r-1]$ is omitted if $i_r = j_r$.
\begin{proposition}[\cite{BZ}]\label{PROP7}
The non-zero quantum minors form a subset of the dual cano\-ni\-cal basis.
More precisely, $\Delta(I,J) = G^*(\m(I,J)).$
\end{proposition} 
In view of Conjecture~\ref{CONJ2'}, it is natural to look for
a necessary and sufficient condition for two quantum minors
to quasi-commute.
This problem was solved in \cite{LZ} in the case of 
quantum flag minors.
We call {\em quantum flag minors} the elements $\Delta(I,J)$ for
which $I$ consists of consecutive integers :
$I=\{i+1,i+2, \ldots ,\,i+k\}=[i+1,i+k]$.
By (\ref{EQQM}), for such a minor and for all $s\ge 0$ we have  
$
\Delta(I,J)=\Delta([i-s,i+k]\,,\,[i-s,i]\cup J).
$
Thus, letting $s$ tend to $\infty$ we get
$\Delta(I,J)=\Delta(\Z_{\le i+k}\,,\,\Z_{\le i}\cup J)$.
Thus a quantum flag minor may be labelled by a unique subset of $\Z$
of the form $\J=\Z_{\le i}\cup J$, where $J$ is a finite subset of
$\Z_{> i}$.
It will be convenient to write $\<\J\>$ in place of $\Delta(I,J)$.

Let $A,\,B$ be finite disjoint subsets of $\Z$ of cardinality $|A|$ and $|B|$. 
We write $A\prec B$ if $a<b$
for any $a\in A,\,b\in B$ (in particular if $A$ or $B$ or both are empty).
More generally, we write $A\Join B$ if one of the following holds:
\begin{quote}
either $|A|\le |B|$ and $A=A'\cup A''$ with $A'\prec B \prec A''$,

or $|B|\le |A|$ and $B=B'\cup B''$ with $B'\prec A \prec B''$.
\end{quote}
Let now $\I= \Z_{\le i}\cup I$ and $\J=\Z_{\le j}\cup J$ be two
infinite subsets as above. Set 
$I':=\I\setminus \J$ and $J':=\J\setminus \I$.
Then $I'$ and $J'$ are finite disjoint subsets. 
We say that $\I$ and $\J$ are {\em separated} if $I'\Join J'$,
and {\em strongly separated} if $I'\prec J'$ or $J'\prec I'$.
We can now state
\begin{proposition}[\cite{LZ}]\label{PROP8}
The quantum flag minors $\<\I\>$ and $\<\J\>$ quasi-commute if and only if
the subsets $\I$ and $\J$ are separated.
\end{proposition}
\begin{example}{\rm Let $a\in\Z$, $\I = \Z_{\le -1} \cup \{1,5\}$ and
$\J(a) = \Z_{\le a} \cup \{a+2,a+4\}$.
Then $\I$ and $\J$ are not separated for $a=\pm 2$, are
separated for $a=0$, and strongly separated in all other cases.
Hence the quantum flag minors $\<\I\>$ and $\<\J(a)\>$ quasi-commute
for all $a \not = \pm 2$.
}
\finex
\end{example}


\section{The Berenstein-Zelevinsky conjecture for flag minors}
\label{SECT6}
\subsection{} \label{SECT6.1}
Consider a product $\pi=\<\I_1\>\cdots \<\I_r\>$ of quantum flag minors.
Recall that the sets $\I_k$ are of the form 
$\I_k = \Z_{\le a_k} \cup I_k$
for some finite set $I_k$ contained in $\Z_{\ge a_k}$.
This decomposition is  not unique, and one could
replace $a_k$ by any $b\le a_k$ and write
$\I_k = \Z_{\le b} \cup ([b+1,a_k] \cup I_k)$.
Hence, taking $b=\min\{a_1,\ldots ,a_r\}$, we can assume
that all $a_k$'s are equal. Moreover, since the translations
of indices $T_{ij} \mapsto T_{i+n,j+n}$ extend to algebra automorphisms
of $A_v$ for all $n\in\Z$, there is no loss of generality
in assuming that $a_k =0$ for all $k$.
So from now on, we will suppose that 
$\I_k = \Z_{\le 0} \cup I_k$
with 
\[
I_k=\{i_1^{(k)}<i_2^{(k)}<\cdots <i_{n_k}^{(k)}\}\subset\Z_{>0}\,. 
\]
The corresponding multi-segments will be denoted by
\[
\m_k = \m([1,n_k],I_k) = \sum_{j=1}^{n_k}\ [j,i_j^{(k)}-1]\,. 
\]
To the product $\pi=\<\I_1\>\cdots \<\I_r\>$ we associate
the multi-segment $\m_\pi=\sum_{k=1}^r \m_k$ and the integer
\[
b_\pi = \sum_{1\le k<l\le r} b(\m_k,\m_l)
\]
where $b(\m,\n)$ was defined in \S \ref{SECT211}. 
We also define the word
\[
w_\pi = i_{n_r}^{(r)} \cdots i_2^{(r)}i_1^{(r)}
i_{n_{r-1}}^{(r-1)} \cdots i_2^{(r-1)}i_1^{(r-1)}
\cdots
i_{n_1}^{(1)} \cdots i_2^{(1)}i_1^{(1)}
\]
obtained by reading the sets $I_r,\ldots,I_1$ successively in
decreasing order.
Recall that the Robinson-Schensted correspondence maps a word $w$
to a pair $(P(w),Q(w))$ of Young tableaux of the same shape $\lambda$
(see for example \cite{Fu}). 
We denote by $\mu_\pi$ the partition conjugate to the shape
of $P(w_\pi)$.
Finally, we introduce the multi-segment
$\n_\pi = \sum_{ i < j} n_{ij}\,[i,j]$
where $n_{ij}$ is the number of letters $j{+}1$ on the $i$th row
of the Young tableau $P(w_\pi)$.
\begin{proposition}\label{TH25}
Suppose that $\mu_\pi$ is equal to the non-increasing reordering 
of $(n_1,\ldots ,n_r)$.
Then 
$
v^{d_\pi}\,\pi \equiv G^*(\n_\pi) \mod v\L^*\,,
$
for some integer $d_\pi$.
\end{proposition}
The proof of Proposition~\ref{TH25} will be given in Section~\ref{SECT7}.
\begin{example}{\rm Let $r=2$ and $I_1=\{2,3,5\},\ I_2=\{1,4\}$.
Then 
\[
\<\I_1\> = G^*([1] + [2]+[3, 4]) \, ,
\qquad \<\I_2\> = G^*( [2, 3]). 
\]
We have 
$w_\pi = 4\,1\,5\,3\,2$ and 
\[
P(w_\pi) =
\matrix{\young{4\cr 3 & 5 \cr  1 & 2\cr}}\,.
\]
(Here and in what follows we choose the French orientation for
drawing Young tableaux.)
Thus, $\mu_\pi=(3,2)=(n_1,n_2)$,  and $\pi$
satisfies the hypothesis of Proposition~\ref{TH25}.
We have $\n_\pi = [1]+[2]+ [3]+ [2, 4]$,
and
\[
v\,\pi = G^*([1]+[2]+[3]+[2,4]) + v\,G^*([1]+[2]+[2,3]+[3,4])\,,
\]
which shows that $v\,\pi \equiv G^*(\n_\pi) \mod v\L^*$.
}
\finex
\end{example} 
A word $w_\pi$ satisfying the hypothesis of Proposition~\ref{TH25} 
is a {\em frank word}, as defined by Lascoux and Sch\"utzenberger \cite{LS}.
The combinatorics of frank words has already occured in several
interesting problems \cite{LS}, \cite{FL}.
Proposition~\ref{TH25} easily yields the following 
\begin{theorem}\label{corol27}
Suppose that the sets $\I_1, \ldots , \I_r$ are pairwise strongly
separated. Then 
\[
v^{b_\pi} \pi = v^{b_\pi} \, \<\I_1\>\cdots \<\I_r\> = G^*(\m_\pi)
\]
belongs to the dual canonical basis.
\end{theorem}
\proof
If two sets $I$ and $J$ are strongly separated then
either 
$I\setminus J\prec J\setminus I$ or 
$J\setminus I\prec I\setminus J$. 
Since the quantum minors $\<\I_k\>$ pairwise quasi-commute, we may 
assume that their indexing is so chosen that 
for $k<l$ we have 
$I_l\setminus I_k\prec I_k\setminus I_l$.
In this case, one checks easily that the $s$th row of the
tableau $P(w_\pi)$ is exactly the reordering of the letters
\[
i_s^{(a_1)},\,i_s^{(a_2)},\,\ldots ,\, i_s^{(a_t)},
\]
where $\{a_1,\ldots ,a_t\}$ is the subset of $[1,r]$
consisting of the integers $k$ for which $n_k \ge s$.
Hence, $\mu_\pi$ is the reordering of $(n_1,\ldots , n_r)$.
It is also clear from this description of $P(w_\pi)$
that $\m_\pi = \n_\pi$ in this case.
By Proposition~\ref{TH25}, it follows that 
\[
v^{d_\pi} \pi \equiv E^*(\m_\pi) \mod v\L^*\,
\]
for some $d_\pi \in \Z$. Now it follows from
(\ref{HIT}) that the expansion of $\pi$ on $\{E^*(\m)\}$
contains $E^*(\m_\pi)$ with coefficient $v^{-b_\pi}$.
Hence we also have $d_\pi = b_\pi$.
Therefore for $r=2$ we have verified Conjecture~\ref{CONJ2'},
and Theorem~\ref{corol27} is proved in this case. 
Otherwise, putting $\pi' = \<\I_1\>\cdots \<\I_{r-1}\>$ we may
suppose by induction on $r$ that 
\[
v^{b_{\pi'}}\,\pi' = G^*(\m_{\pi'})\,.
\]
Then $G^*(\m_{\pi'})$ and $\<\I_r\>=G^*(\m_r)$ quasi-commute and
satisfy
\[
v^{b(\m_{\pi'},\m_r)}\,G^*(\m_{\pi'})\,G^*(\m_r) \equiv
E^*(\m_{\pi'}+\m_r) \mod v\L^*\,.
\]
Hence $G^*(\m_{\pi'})$ and $G^*(\m_r)$ 
verify Conjecture~\ref{CONJ2'}, and the result is proved. 
\cqfd

\subsection{}\label{SECT46N}
In the case of a product of $r=2$ quantum flag minors, we can 
drop the hypothesis of {\em strong} separation in Theorem~\ref{corol27}
and deduce from Proposition~\ref{TH25} the following
\begin{theorem}\label{corol28} 
Let $\pi = \<\I\>\<\J\>$ be a product of two 
quantum flag minors. Then $\pi$ 
belongs to the dual canonical basis up to a power of $v$
if and only if $\<\I\>$ and $\<\J\>$ quasi-commute. 
\end{theorem}
The proof of Theorem~\ref{corol28} will be given in
Section~\ref{SECT7}.

\subsection{} Finally, using an argument of \cite{KMT,MT} about
finite-dimensional representations of $U_v(\widehat{\Sl}_N)$
as well as Theorem~\ref{corol28}, we can improve on Theorem~\ref{corol27}
and obtain the main result.
\begin{theorem}\label{MainTH} 
Let $\pi=\<\I_1\>\cdots \<\I_r\>$ be a product of any number $r$ of
quantum flag minors. Then $\pi$
belongs to the dual canonical basis up to a power of $v$ if
and only if the minors $\<\I_k\>$ pairwise quasi-commute. 
\end{theorem}
Note that if $\pi$ belongs to $\B_v$ up to a power of $v$,
this power has to be $v^{b_\pi}$.
The proof of Theorem~\ref{MainTH} will be given in
Section~\ref{SECT8}.


\section{Proofs of Proposition~\ref{TH25} and Theorem~\ref{corol28}}
\label{SECT7}

\subsection{} We retain the notation of \ref{SECT6.1}. 
Let $N=\max(I_1 \cup \cdots \cup I_r)$.
The quantum minors $\<\I_1\>, \ldots , \<\I_r\>$ belong to the subalgebra
of $U_v(\nn_\infty)$ generated by $E_1, \ldots , E_{N-1}$, which is
isomorphic to $U_v^+(\Sl_N)$.
So Proposition~\ref{TH25} can be regarded as a statement about $U_v(\g)$,
where $\g = \Sl_N$.
We  denote by 
$E_i, F_i, K_i\ (i = 1,\ldots , N-1)$
the standard generators of $U_v(\g)$.
The subalgebras generated by $E_i\ (i = 1,\ldots , N-1)$
(\resp by $E_i, K_i\ (i = 1,\ldots , N-1)$)
are denoted by
$U_v(\nn)$ (\resp $U_v(\bb)$). 
The fundamental weights of $\g$ are denoted by 
$\Lambda_1, \ldots , \Lambda_{N-1}$.

\subsection{} We recall some general properties of the canonical bases
of $U_v(\nn)$ and of the finite-dimensional irreducible $U_v(\g)$-modules.
Our aim is to express the multiplication of vectors of the dual
canonical basis of $U_v(\nn)$ in terms of the tensor products of vectors
of the dual canonical bases of the irreducible modules.
The precise relation will be given by \ref{CONCL} (\ref{EQ19}) below.

\subsubsection{} Let $\lambda$ be an integral dominant weight of $\g$
and $V(\lambda)$ the irreducible $U_v(\g)$-module with highest
weight $\lambda$.
Let $\M_N$ denote the set of multi-segments supported on $[1,N{-}1]$.
The canonical (or lower global) basis of $U_v(\nn)$ may be identified
with the subset of the canonical basis
of $U_v(\nn_\infty)$ consisting of the vectors $G(\m)$ with $\m\in\M_N$.
Let $u_\lambda$ be a fixed lowest weight vector of $V(\lambda)$ (hence
$u_\lambda$ has weight $w_0\lambda$ where $w_0$ is the longest element
of the Weyl group of $\g$).
It is known that the image of $\{G(\m) \ | \ \m \in \M_N\}$ under the
map
\[ 
\pi_\lambda : U_v(\nn) \longrightarrow V(\lambda),
\quad x \longrightarrow xu_\lambda
\]
is the union of a basis of $V(\lambda)$ with the set $\{0\}$.
This basis is called the canonical (or lower global) basis of
$V(\lambda)$ and is naturally labelled by the set $\tab(\lambda)$
of Young tableaux of shape $\lambda$ over $[1,N]$. 
(Here we identify as usual the dominant integral weight $\lambda$
with a partition of length $\le N-1$.)
This basis will be denoted by $\{G(t) \ | \ t\in\tab(\lambda)\}$.

\subsubsection{} \label{DUAL} Given a left $U_v(\g)$-module $M$, we endow
$M^* = \hom(M,\Q(v))$ with a left $U_v(\g)$-action by setting
\[
(x\phi)(m) = \phi(x^\natural m)\, \qquad x \in U_v(\g),\ m\in M,\ \phi
\in M^*,
\]
where $x \mapsto x^\natural$ is the anti-automorphism of $U_v(\g)$ defined by
\[
E_i^\natural = F_i, \quad F_i^\natural = E_i, \quad K_i^\natural = K_i\,.
\]
It is known that $V(\lambda)^*$ is isomorphic to $V(\lambda)$ as 
a $U_v(\g)$-module. 
In other words, there is a non-degenerate scalar product on 
$V(\lambda)$ satisfying 
\[
(u_\lambda,u_\lambda) = 1,\qquad
(xu,v) = (u,\sigma(x)v), \quad x\in U_v(\g), \ u,v \in V(\lambda),
\]
and $V(\lambda)^*$ can be identified to $V(\lambda)$ via the map
$V(\lambda) \longrightarrow V(\lambda)^*$, 
$x \longrightarrow ( x\,, \cdot )$. 

The basis $\B_v(\lambda) =\{G^*(t) \ | \ t\in \tab(\lambda) \}$
dual to $\{G(t)\}$ is called the dual canonical (or upper global)
basis of $V(\lambda)^*$ (or $V(\lambda)$).

\subsubsection{} \label{DUALN} Similarly, Kashiwara's scalar product 
on $U_v(\nn)$
(defined as in \ref{SECT24}) allows to identify $U_v(\nn)$ to its
graded dual $U_v(\nn)^*$, as vector spaces.
By composing these maps we obtain an embedding 
\[
\psi_\lambda : V(\lambda) \simeq V(\lambda)^*
\stackrel{\pi_\lambda^*}\longrightarrow
U_v(\nn)^* \simeq U_v(\nn)
\]
satisfying $\psi_\lambda(u_\lambda) =1$ and 
$\psi_\lambda(\B_v(\lambda)) \subset \B_v$, where, as in \ref{SECT24},
$\B_v = \{G^*(\m) \ | \ \m \in \M_N\}$ is the basis dual to the 
canonical basis $\{G(\m)\ |\ \m \in \M_N\}$.

More precisely, write $\lambda = \sum_{i=1}^{N-1} a_i \Lambda_i$ and set 
$\tilde{\lambda} = \sum_{i=1}^{N-1} a_i \Lambda_{N-i}$.
If $t\in\tab(\lambda)$ has columns $C_1,\ldots ,C_m$, where 
$m=\sum_i a_i$, we define $\tilde{t}\in\tab(\tilde{\lambda})$
as the Young tableau with columns $D_1, \ldots , D_m$, where 
$D_i = [1,N] \setminus C_{m+1-i}$. (Here, columns of Young
tableaux are identified to subsets of $[1,N]$ in a natural way).
Finally, we introduce the multi-segment 
$\m(t)= \sum_{i<j} m_{ij}(t) \, [i,j]$,
where $m_{ij}(t)$ is the number of letters $j+1$ in the $i$th
row of $\tilde{t}$.
Then $\psi_\lambda(G^*(t)) = G^*(\m(t))$.

\subsubsection{}\label{724}
Let $\lambda$ and $\mu$ be two dominant weights. The tensor product
$V(\lambda) \otimes V(\mu)$ is a $U_v(\g)$-module via the
comultiplication
\begin{eqnarray}
  \Delta E_i &=& E_i \otimes 1 + K_i^{-1} \otimes E_i,\label{DELTAE} \\
  \Delta F_i &=& F_i \otimes K_i + 1 \otimes F_i, \\
  \Delta K_i &=& K_i \otimes K_i,  
\end{eqnarray}
and we have the following commutative diagram
\[
\begin{array}{ccc}
U_v(\g) & \stackrel{\Delta}\longrightarrow & U_v(\g) \otimes U_v(\g)
\\ [2mm]
p_{\lambda+\mu}\downarrow & & \downarrow p_\lambda \otimes p_\mu \\
V(\lambda+\mu) & \stackrel{i_{\lambda,\mu}}\longrightarrow &
V(\lambda)\otimes V(\mu)
\end{array}
\]
where $p_\lambda : U_v(\g) \longrightarrow V(\lambda), \
x \longrightarrow xu_\lambda$, and 
$i_{\lambda,\mu}$ is the homomorphism of $U_v(\g)$-modules
mapping $u_{\lambda+\mu}$ to $u_\lambda\otimes u_\mu$.
Dualizing, we obtain
\[
\begin{array}{ccc}
U_v(\g)^*\otimes U_v(\g)^* & \stackrel{m}\longrightarrow & U_v(\g)^* 
\\ [2mm]
p_{\lambda}^*\otimes p_{\mu}^*\uparrow & & \uparrow p_{\lambda+\mu}^* \\
V(\lambda)^*\otimes V(\mu)^* & \stackrel{i_{\lambda,\mu}^*}\longrightarrow &
V(\lambda+\mu)^*
\end{array}
\]
where $m$ is the multiplication dual to $\Delta$, and 
$V(\lambda)^*\otimes V(\mu)^* \simeq (V(\lambda)\otimes V(\mu))^*$
is endowed with the action induced by the comultiplication
\begin{eqnarray}
  \Delta' E_i &=& E_i \otimes K_i + 1 \otimes E_i, \label{EQD'1} \\
  \Delta' F_i &=& F_i \otimes 1 + K_i^{-1} \otimes F_i, \label{EQD'2}\\
  \Delta' K_i &=& K_i \otimes K_i,  \label{EQD'3}
\end{eqnarray}
(this comes from \ref{DUAL}).

\subsubsection{} In contrast to the case $v=1$, $U_v(\nn)$ is not a 
Hopf subalgebra of $U_v(\g)$ : it is not stable under $\Delta$.
Nevertheless, one can define a multiplication on $U_v(\nn)^*$ by
embedding it into
$U_v(\bb)^*$ via $\phi \longrightarrow \tilde{\phi}$, where
for $\phi\in U_v(\nn)^*$, 
\[
\tilde\phi(xK_i) = \phi(x), \qquad x\in U_v(\nn),\ i=1,\ldots , N-1.
\]
It is easy to see that the image of $U_v(\nn)^*$ is a subalgebra of
$U_v(\bb)^*$, hence $U_v(\nn)^*$ inherits a multiplication.
It is also known that the vector space isomorphism 
$U_v(\nn) \longrightarrow U_v(\nn)^*$ of \ref{DUALN} becomes in this 
way an algebra isomorphism. 

\subsubsection{} The image 
of $V(\lambda)^*$ under $p_\lambda^*$ can be regarded as a subset of
$U_v(b)^*$.
More precisely, let $\phi \in V(\lambda)^*$ and set
$\overline{\phi}=\pi_\lambda^*(\phi)\in U_v(\nn)^*$,
$\tilde{\phi}=p_\lambda^*(\phi)\in U_v(\bb)^*$.
For $x\in U_v(\nn)$, we have 
\[
\tilde{\phi}(xK_i) = \phi(xK_iu_\lambda) 
=v^{(\alpha_i,w_0\lambda)}\,\phi(xu_\lambda)
=v^{(\alpha_i,w_0\lambda)}\,\overline{\phi}(x).
\]
Let $Q^+ = \bigoplus_{i=1}^{N-1} \N \alpha_i$ be the monoid
generated by the simple roots of $\g$.
We say that $\overline{\phi}$ (or $\tilde{\phi}$) has weight 
$\alpha \in Q^+$ if $\overline{\phi}$ vanishes on all homogeneous
elements of $U_v(\n)$ of weight $\not = \alpha$.
\begin{lemma}\label{LEMPRODUIT}
Let $\phi \in V(\lambda)^*$, $\psi \in V(\mu)^*$ and suppose that
$\overline{\psi}$ has weight $\beta$. We have
\[
\left(\tilde{\phi}\cdot\tilde{\psi}\right)(xK_i) 
=
v^{(\alpha_i\,,\,w_0(\lambda + \mu)) - (\beta\,,\,w_0\lambda)}\,
\left(\overline{\phi}\cdot\overline{\psi}\right)(x),
\qquad x\in U_v(\nn),\ i=1,\ldots N-1,
\]
where $\tilde{\phi}\cdot\tilde{\psi}$ denotes multiplication in
$U_v(\bb)^*$ and $\overline{\phi}\cdot\overline{\psi}$ multiplication
in $U_v(\nn)^*$. 
\end{lemma}
\proof
Let $x\in U_v(\nn)$ and put 
$\Delta x = \sum_{(x)} x_{(1)} \otimes x_{(2)}$, where 
$x_{(1)}\in U_v(\bb)$ and $x_{(2)} \in U_v(\nn)$.
In view of (\ref{DELTAE}), one can moreover assume that $x_{(2)}$ is
homogeneous of weight $\sum_{i=1}^{N-1} a_i\alpha_i$, and 
$x_{(1)}=x'_{(1)}\,\prod_{i=1}^{N-1}K_i^{-a_i}$
with $x'_{(1)}\in U_v(\nn)$.
Then
\begin{eqnarray*}
 \left(\tilde{\phi}\cdot\tilde{\psi}\right)(xK_i) & = &
\tilde{\phi}\otimes\tilde{\psi}(\Delta(xK_i)) \\
&=&
\sum_{(x)}\tilde{\phi}(x_{(1)}K_i)\,\tilde{\psi}(x_{(2)}K_i)\\
&=&
v^{(\alpha_i\,,\,w_0(\lambda+\mu))}\,
\sum_{(x)}\tilde{\phi}\left(x'_{(1)}\prod
K_i^{-a_i}\right)\tilde{\psi}(x_{(2)})\\
&=&
v^{(\alpha_i\,,\,w_0(\lambda+\mu))}\,
v^{(-\beta\,,\,w_0\lambda)}\,
\sum_{\wt(x_{(2)})=\beta}
\overline{\phi}(x'_{(1)})\,\overline{\psi}(x_{(2)})\\
&=&
v^{(\alpha_i\,,\,w_0(\lambda+\mu))}\,
v^{(-\beta\,,\,w_0\lambda)}\,
\left(\overline{\phi}\cdot\overline{\psi}\right)(x).
\end{eqnarray*}
\cqfd

\subsubsection{} \label{CONCL}
It follows from \ref{724} and Lemma~\ref{LEMPRODUIT} that if $q_{\lambda,\mu}$
denotes the homomorphism of $U_v(g)$-modules
\[
q_{\lambda,\mu} : V(\lambda)\otimes V(\mu) \simeq
V(\lambda)^*\otimes V(\mu)^*
\stackrel{i_{\lambda,\mu}^*}\longrightarrow
V(\lambda+\mu)^* \simeq V(\lambda + \mu)\,,
\]
then for $y\in V(\lambda)$ and $z\in V(\mu)$,
\[
\psi_{\lambda+\mu}\left(q_{\lambda,\mu}(y\otimes z)\right)
=
v^{(w_0\lambda \,,\,-\wt(z) +w_0\mu )}\,
\psi_\lambda(y) \cdot \psi_\mu(z)\, .
\]
In particular, taking $t\in \tab(\lambda)$, $s\in\tab(\mu)$ and 
$y=G^*(t)$, $z=G^*(s)$, by \ref{DUALN} we have $\psi_\lambda(y) = G^*(\m)$,
$\psi_\mu(z) = G^*(\n)$ for certain $\m,\,\n \in \M_N$, and
\begin{equation}\label{EQ19}
\psi_{\lambda+\mu}\left(q_{\lambda,\mu}(G^*(t)\otimes G^*(s))\right)
=
v^{-(w_0\lambda \,,\,\wt(\n))}\,
G^*(\m)\, G^*(\n)\, .  
\end{equation}
 
\subsection{} We now come to the proof of Proposition~\ref{TH25}. 
\subsubsection{} Let $(L(\lambda),B(\lambda))$ be the upper crystal
basis of $V(\lambda)$ \cite{Ka1,Ka2,Ka3}.
We recall that the upper crystal lattice $L(\lambda)$ is the 
${\bf A}$-span of $\{G^*(t) \ | \ t\in \tab(\lambda) \}$,
and that $B(\lambda)$ is the $\Q$-basis 
$\{G^*(t) \mod vL(\lambda) \ | \ t\in \tab(\lambda) \}$
of the $\Q$-vector space $L(\lambda)/vL(\lambda)$.
The elements of $B(\lambda)$ can be seen as combinatorial
labels for the vectors of the (dual) canonical basis.
We shall identify them with Young tableaux by writing $t$
in place of $G^*(t) \mod vL(\lambda)$.

\subsubsection{} \label{732} The comultiplication $\Delta'$ of 
(\ref{EQD'1}), (\ref{EQD'2}), (\ref{EQD'3})
is compatible with upper crystal bases \cite{Ka1},
therefore 
$q_{\lambda,\mu}(L(\lambda)\otimes L(\mu))
= L(\lambda + \mu)$.
Moreover, if $t\otimes s$ belongs to the unique connected component
of the crystal graph of $V(\lambda)\otimes V(\mu)$ isomorphic
to $B(\lambda + \mu)$, then 
\[
q_{\lambda,\mu}(G^*(t)\otimes G^*(s)) \equiv G^*(t\cdot s)
\quad \mod vL(\lambda + \mu)
\]
for some element $t\cdot s$ of $B(\lambda + \mu) = \tab(\lambda + \mu)$.

\subsubsection{} \label{733}
For $t\in\tab(\lambda)$, let $w(t)$ be the word obtained by
reading the columns of $t$ in decreasing order and from left to right.
Thus, if 
\[
t =\matrix{\young{6\cr 4 & 5 \cr  2 & 3\cr 1 & 1\cr}}
\]
then $w(t) = 6\,4\,2\,1\,5\,3\,1$.
Recall the Robinson-Schensted map $w \mapsto (P(w),Q(w))$ \ref{SECT6.1}.
It is known (see \cite{LT}) that for 
$t\in\tab(\lambda),s\in\tab(\mu)$,
the product $t\otimes s$ belongs to the connected component
of the crystal of $V(\lambda)\otimes V(\mu)$ of type 
$B(\lambda +\mu)$ if and only if the concatenation $w(t)\,w(s)$
of the words $w(t)$ and $w(s)$ is mapped under the Robinson-Schensted
map to a Young tableau $P(w(t)\,w(s))$ of shape $\lambda + \mu$. 
In this case $t\cdot s = P(w(t)\,w(s))$.

\subsubsection{} Using \ref{CONCL} we can now prove
\begin{proposition}\label{PROP34}
Let $\m, \ \n \in \M_N$.
Suppose there exist $t\in\tab(\lambda),\ s\in\tab(\mu)$
such that $\m = \m(t),\ \n = \m(s)$ and $t\otimes s$
belongs to the connected component of $B(\lambda)\otimes B(\mu)$
of type $B(\lambda + \mu)$.
Then
\[
v^{-(w_0\lambda \,,\,\wt(\n))}\,G^*(\m)\,G^*(\n)
\equiv
G^*(\p) \quad \mod v\L^*\,,
\]
where $\p = \m(t\cdot s)$.
\end{proposition}
\proof
By (\ref{EQ19}) we know that
\[
v^{-(w_0\lambda \,,\,\wt(\n))}\,G^*(\m)\,G^*(\n)
= \psi_{\lambda + \mu}\left(q_{\lambda,\mu}(G^*(t)\otimes G^*(s))\right)\,.
\]
By \ref{732}, 
$
q_{\lambda,\mu}(G^*(t)\otimes G^*(s)) \equiv G^*(t\cdot s)
\mod vL(\lambda + \mu) \,.
$
Since for $u\in\tab(\lambda + \mu)$, by \ref{DUALN},
$\psi_{\lambda + \mu}(G^*(u)) = G^*(\m(u))$, we deduce that
\[
\psi_{\lambda + \mu}\left(q_{\lambda,\mu}(G^*(t)\otimes G^*(s))\right)
\equiv
G^*(\m(t\cdot s)) \quad \mod vL^* \,,
\]
where $L^*$ denotes the ${\bf A}$-lattice spanned by the dual
canonical basis of $U_v(\nn)$. Hence we have proved
\[
v^{-(w_0\lambda \,,\,\wt(\n))}\,G^*(\m)\,G^*(\n)
\equiv
G^*(\p) \quad \mod vL^*\,.
\]
On the other hand we know that $G^*(\m)\,G^*(\n)$ is a 
$\Z[v,v^{-1}]$-linear combination of vectors of the
dual canonical basis. Therefore the congruence is indeed
modulo $v\L^*$. (Recall that $\L^*$ is the $\Z[v]$-lattice spanned
by the dual canonical basis.)
\cqfd
\subsubsection{}
Suppose now that $G^*(\m)$ and $G^*(\n)$ are quantum flag minors,
namely, that $\m = \m_1,\ \n=\m_2$ are associated with subsets 
$I_1,\ I_2$ of cardinality $n_1,\ n_2$ in the notation of \ref{SECT6.1}.
Then $G^*(\m_1) = \psi_{\Lambda_{N-n_1}}(G^*(t_1))$,
$G^*(\m_2) = \psi_{\Lambda_{N-n_2}}(G^*(t_2))$, where $t_1,\ t_2$
are the Young tableaux of column shape associated with the sets
$[1,N] \setminus I_1,\ [1,N] \setminus I_2$.
Applying Proposition~\ref{PROP34} in this case, we see that we 
obtain precisely Proposition~\ref{TH25} for the product of two
quantum flag minors.

\subsubsection{} The extension to a product of $r>2$ minors is 
straightforward.
First, \ref{CONCL} generalizes easily as follows.
If $t_i\in\tab(\lambda^{(i)})$ and 
$\psi_{\lambda^{(i)}}(G^*(t_i)) = G^*(\m_i)\ 
(i=1,\ldots ,r)$, then, denoting by 
$q_{\lambda^{(1)},\ldots ,\lambda^{(r)}}$
the projection
$V(\lambda^{(1)})\otimes \cdots \otimes V(\lambda^{(r)})
\longrightarrow
V(\lambda^{(1)}+\cdots + \lambda^{(r)})$
we have
\[
\psi_{\lambda^{(1)}+\cdots +\lambda^{(r)}}
(q_{\lambda^{(1)},\ldots ,\lambda^{(r)}}
(G^*(t_1)\otimes \cdots \otimes G^*(t_r)))
=
v^{-\sum_{i<j}(w_0\lambda^{(i)} ,\wt(\m_j))}
G^*(\m_1)\cdots G^*(\m_r) .  
\]
Secondly, \ref{732} and \ref{733} generalize in the obvious way,
namely, if $t_1\otimes \cdots \otimes t_r$ belongs to the
connected component of 
$B(\lambda^{(1)}) \otimes \cdots \otimes B(\lambda^{(r)})$
of type $B(\lambda^{(1)}+\cdots \lambda^{(r)})$, then
\[
q_{\lambda^{(1)},\ldots ,\lambda^{(r)}}
(G^*(t_1)\otimes \cdots \otimes G^*(t_r))
\equiv
G^*(t_1 \cdots  t_r)
\quad \mod vL(\lambda^{(1)}+\cdots +\lambda^{(r)})
\]
where
$
t_1 \cdots  t_r
=
P(w(t_1) \cdots  w(t_r)).
$
From these facts one obtains easily the generalization of 
Proposition~\ref{PROP34} to $r$ factors, and by specializing
to quantum flag minors one gets the general case of
Proposition~\ref{TH25}.
This finishes the proof of Proposition~\ref{TH25}.

\subsection{} In this section we prove Theorem~\ref{corol28}.
\subsubsection{}\label{IDENTIF}
Let us recall the realization of the finite-dimensional irreducible
$U_v(\g)$-modules $V(\lambda)$ as subspaces of
the deformation $A_v[\F]$ of the coordinate 
ring of the flag variety \cite{TT,LT,LZ}.
The algebra $A_v[\F]$ is generated over $\Q(v)$ by $2^N - 1$ generators
$[J]$ labelled by the non-empty subsets of $[1,N]$ submitted to the
relations $R_1$ and $R_2$ below.
Given two subsets $A$ and $B$ of $[1,N]$, write 
$\inv(A,B)$ for the number of pairs $(a,b)\in A\times B$ with $a>b$.
\begin{quote}
$R_1(I,J)$: for all subsets $I$ and $J$ such that $|I| \le |J|$, we
have
\[
[I][J] = \sum_M (-v)^{\inv(J\setminus M,M)-\inv(I,M)}\,[I\cup M][J \setminus
M]\,,
\]
where the sum is over all $M\subset J\setminus I$ with $|M| =
|J|-|I|$;

$R_2(I,J)$:  for all subsets $I$ and $J$ such that $|I|-1 \ge |J|+1$, we
have  
\[
\sum_{i\in I\setminus J} 
(-v)^{\inv(\{i\},I\setminus \{i\}) - \inv(\{i\},J)}\,
[I\setminus \{i\}][J\cup \{i\}] = 0\,.
\]
\end{quote}
To each subset $I$ of cardinality $r \le N$ we associate the vector
$G^*(I)$ of the fundamental representation $V(\Lambda_r)$.
Let $(i_1,\ldots , i_r) \in [1,N-1]^r$ and set 
$\lambda = \sum_{k=1}^r\Lambda_{i_k}$.
(Here, as in \ref{DUALN}, we identify $I\subset [1,N]$ to
a column Young tableau.)
Define a linear map $Q_{\Lambda_{i_1},\ldots ,\Lambda_{i_r}}$
from $V(\Lambda_{i_1})\otimes \cdots \otimes V(\Lambda_{i_r})$
to $A_v[\F]$ by
\[
Q_{\Lambda_{i_1},\ldots ,\Lambda_{i_r}}
\left(G^*(I_1)\otimes \cdots \otimes G^*(I_r)\right)
=
[I_1]\cdots [I_r]\,.
\]
Then $\ker(Q_{\Lambda_{i_1},\ldots ,\Lambda_{i_r}})
= \ker(q_{\Lambda_{i_1},\ldots ,\Lambda_{i_r}})$,
hence $Q_{\Lambda_{i_1},\ldots ,\Lambda_{i_r}}$ induces an embedding
of $V(\lambda)$ in $A_v[\F]$ which maps 
$q_{\Lambda_{i_1},\ldots ,\Lambda_{i_r}}
(G^*(I_1)\otimes \cdots \otimes G^*(I_r))$
to $[I_1]\cdots [I_r]$.
In other words the linear relations in $V(\lambda)$ between the vectors
$q_{\Lambda_{i_1},\ldots ,\Lambda_{i_r}}
(G^*(I_1)\otimes \cdots \otimes G^*(I_r))$
are exactly the same as those in $A_v[\F]$ between the vectors
$[I_1]\cdots [I_r]$.

\subsubsection{} 
Let $I$ and $J$ be two separated subsets of $[1,N]$ 
and write 
$I_0 = I\setminus (I\cap J),\ 
J_0 = J\setminus (I\cap J)$. 
Up to a possible exchange of $I$ and $J$,
we can assume that $|I|\le |J|$ and $I_0 = I'\cup I''$
with $I' \prec J_0 \prec I''$.
\begin{proposition}\label{PROP20F}
The following relation holds in $A_v[\F]$:
\begin{equation}\label{EQFM}
[J][I] = v^{|I'|}\,\sum_M (-v)^{\inv(M,J_0\setminus M)}\,
[I'\cup (J\setminus M)]\,
[M\cup(I\setminus I')] \,,
\end{equation}
where the sum is over all $M\subset J_0$ with $|M|=|I'|$.
\end{proposition}
\proof
First, because of \cite{LZ} Lemma~2.3, it is enough to prove the
proposition in the case $I\cap J=\emptyset$, so we can assume
that $I_0 = I$ and $J_0 = J$.
As a second reduction, we remark that the subset $I''$ is contained
in all right factors of the identity, and that we have $I''\succ I'$
and $I''\succ M$.
Thus, by using the $v$-analogue of Laplace's expansion 
(\cite{TT} Proposition~2.10) we may erase $I''$ in all terms
and we are reduced to prove
\[
[J][I'] = v^{|I'|}\,\sum_M (-v)^{\inv(M,J\setminus M)}\,
[I'\cup (J\setminus M)]\,
[M] \,.
\]
Now since $I'\prec J$ and $|I'|\le|J|$ we have
$[J][I'] = v^{|I'|}\,[I'][J]$ and the identity to be proved becomes
\[
[I'][J] = \sum_M (-v)^{\inv(M,J\setminus M)}\,
[I'\cup (J\setminus M)]\,[M] \,.
\]
Writing $M'= J\setminus M$, we can rewrite it as
\[
[I'][J] = \sum_{M'} (-v)^{\inv(J\setminus M', M')}\,
[I'\cup M']\,[J \setminus M'] \,.
\]
This is exactly $R_1(I',J)$, since $I'\prec J$ implies $\inv(I',M')=0$
for all $M'\subset J$.
\cqfd
 
\begin{example}{\rm In $A_v[\F]$ we have
\[
\begin{array}{rcl}
[3,4,5,6]\,[1,2,7] &=& v^2\,[1,2,5,6]\,[3,4,7] -v^3\,[1,2,4,6]\,[3,5,7]
                  +v^4\,[1,2,3,6]\,[4,5,7] \\[1mm]
&&+\,v^4\,[1,2,4,5]\,[3,6,7]-v^5\,[1,2,3,5]\,[4,6,7]
+v^6[1,2,3,4]\,[5,6,7]\,.
\end{array}
\]
}
\finex
\end{example}

\subsubsection{} Write $|I|=i$ and $|J|=j$. 
Using \ref{IDENTIF} we identify 
$V(\Lambda_i+\Lambda_j)$ 
with the subspace spanned by the vectors $[K][L]$ with $|K|=j$ and $|L|=i$.
We can then interpret (\ref{EQFM}) as a relation between
the vectors 
$q_{\Lambda_j,\Lambda_i}(G^*(K)\otimes G^*(L))$.
Note that all the products $[K][L]$ occuring in the right-hand side
of (\ref{EQFM}) have the property that if we 
arrange the
elements of $K$ and $L$ in increasing order into two columns
shape Young tableaux $t(K)$ and $t(L)$, the 
juxtaposition $t(K)t(L)$ is a semi-standard Young tableau $t(K,L)$.
Therefore, by \ref{732}, \ref{733}, we have 
\[
[K][L] \equiv G^*(t(K,L)) \quad \mod vL(\Lambda_j+\Lambda_i) \,. 
\]
Moreover, the term $[K][L]$ having the lowest power of $v$
is the one in which $M$ consists of the $|I'|$ smaller elements
of $J_0$, and the corresponding tableau $t_{\rm low}$ is obtained
by reordering the rows of the juxtaposition $t(I)t(J)$.
It follows that, using the notation of \ref{DUALN},
$\m(t_{\rm low}) = \m(t(I)) + \m(t(J))$. 
Define
$$
\overline{I} = [1,N] \setminus I,\quad
\overline{J} = [1,N] \setminus J,\quad
\I = \Z_{\le 0} \cup \overline{I},\quad
\J = \Z_{\le 0} \cup \overline{J}.
$$ 
Then, by \ref{DUALN}, 
$$
\psi_{\Lambda_i}([I]) = \<\I\> = G^*(\m(t(I))),\quad
\psi_{\Lambda_j}([J]) = \<\J\> = G^*(\m(t(J)))
$$
and
$$
\psi_{\Lambda_j+\Lambda_i}(G^*(t_{\rm low})) = G^*(\m(t(I))+\m(t(J))).
$$
Hence, applying $\psi_{\Lambda_j+\Lambda_i}$ to both sides of (\ref{EQFM})
we get
\[
v^k\,\<\J\>\<\I\> = v^k\, G^*(\m(t(J)))\,G^*(\m(t(I)))
\equiv
G^*\left(\m(t(I))+\m(t(J))\right) \quad \mod v\L^*\,,
\]
for some integer $k$.
Since we know that $G^*\left(\m(t(I))+\m(t(J))\right)$
occurs in $\<\J\>\<\I\>$ with coefficient 
$v^{-b(\m(t(J)),\m(t(I)))}$, we see that
$k=b(\m(t(J)),\m(t(I)))$.
Note that $I$ and $J$ are separated if and only if $\overline{I}$
and $\overline{J}$ are separated.
Thus, using the reformulation of the Berenstein-Zelevinsky 
conjecture given by Conjecture~\ref{CONJ2'}, we have obtained that
if $\I$ and $\J$ are separated,
$v^{b(\m(t(J)),\m(t(I)))}\,\<\J\>\<\I\>$
belongs to the dual canonical basis $\B_v$.
This finishes the proof of Theorem~\ref{corol28}.


\section{Irreducibility of induction products of evaluation modules} 
\label{SECT8}

\subsection{}
Let $a$ be an integer, 
$\alpha = (\alpha_1,\ldots ,\alpha_r)$ a partition of $m$,
and let $S_\alpha(t^a)$ denote the corresponding evaluation 
module of the affine Hecke algebra $\H_m(t)$.
Recall from \ref{SECT13} that 
$S_\alpha(t^a) \simeq L_{\m(\alpha,a)}$, where
$\m(\alpha,a) = \sum_{i=1}^r [a-i+1,a-i+\alpha_i]$.
Therefore, by Theorem~\ref{THAR}, the class of $S_\alpha(t^a)$ 
in the Grothendieck ring
$R$ is mapped under the isomorphism $\Psi$ to the vector 
$g^*(\m(\alpha,a))$ of $A$. 
By Proposition~\ref{PROP7}, the corresponding element of $A_v$
is
\[
G^*(\m(\alpha,a))=\Delta([a-r+1,a],\{a-r+1+\alpha_r,\ldots ,a+\alpha_1\})\,
\]
that is, a quantum flag minor. 
It follows that we can apply Proposition~\ref{PROP8}, Theorem~\ref{TH10} 
and Theorem~\ref{corol28} and get a completely explicit answer to the 
problem of the irreducibility of the induction product of two
evaluation modules.

Let $\beta=(\beta_1,\ldots ,\beta_s)$ be another partition, and let $b\in\Z$.
We set 
\[
\I:=\Z_{\le a-r} \cup \{a-r+1+\alpha_r,\ldots ,a+\alpha_1\}, \quad
\J:=\Z_{\le b-s} \cup \{b-s+1+\beta_s,\ldots ,b+\beta_1\}.
\]
Then we have $G^*(\m(\alpha,a)) = \<\I\>$ and 
$G^*(\m(\beta,b)) = \<\J\>$ in the notation of Section~\ref{SECT5}.
Hence, we get
\begin{theorem}\label{TH23}
The induction product $S_\alpha(t^a)\odot S_\beta(t^b)$
is irreducible if and only if the sets $\I$ and $\J$ are 
separated.
\end{theorem}

\begin{example}{\rm Take $\alpha = (4,2)$ and $\beta=(2,2,1)$.
We want to determine for which values of $a$ and $b$ the induction
product $S_\alpha(t^a)\odot S_\beta(t^b)$ is irreducible. 
Clearly, this
only depends on the difference $b-a$ and we can take $a=0$.
We have
\[
\I=\Z_{\le -2} \cup \{1,4\},\quad
\J=\Z_{\le b-3} \cup \{b-1,b+1,b+2\}\,.
\]
For $b\le -4$ we have $\J \subset \I$, and for $b\ge 7$, $\I \subset \J$.
Thus if $b \not \in [-3,6]$, $\I$ and $\J$ are trivially separated
and the product is irreducible.
For the remaining values we find that $\I$ and $\J$ are  separated
for $b=0$, strongly separated for $b=2,\,5$, and not separated for
$b=-3,\,-2,\,-1,\,1,\,3,\,4,\,6$.
Hence $S_\alpha(t^a)\odot S_\beta(t^b)$ is simple if and only if
$b-a \not \in \{-3,\,-2,\,-1,\,1,\,3,\,4,\,6\}$.
}
\finex
\end{example}
It is easy to deduce from the definition of separated sets given in
\ref{SECT29} that Theorem~\ref{TH23} is equivalent to the statement of
Theorem~\ref{MAIN1} in the Introduction.
Let us check it.
We first note that the sets $\I\setminus\J$ and $\J\setminus\I$ are
finite and that we have the following relation between their
cardinalities:
$
| \I\setminus\J | = | \J\setminus\I | + a - b.
$ 
Therefore, if $a-b>0$, then $| \I\setminus\J | > | \J\setminus\I |$
and the sets $\I$ and $\J$ are not separated if and only if there
exist $i,j,k$ such that $i,k \in \I\setminus\J$, $j\in  \J\setminus\I$
and $i<j<k$ (see \cite{LZ}, 3.8). The case $a-b<0$ is similar. 
If $a=b$, the two sets $\I\setminus\J$ and $\J\setminus\I$
being of equal cardinality, they are not separated if and only if
there exist $i,j,k,l$ such that  $i,k \in \I\setminus\J$, $j,l\in  \J\setminus\I$
and either $i<j<k<l$ or $j<i<l<k$ (\cite{LZ}, 3.8).

\smallskip
We note the following obvious corollary of Theorem~\ref{MAIN1} 
which will be used in the proof of Theorem~\ref{MAIN2}. 
\begin{corollary} Let $L=S_\alpha(z)$ be an evaluation module
of $\H_m$. Then, $L\odot L$ is irreducible. 
\end{corollary}
\begin{figure}[t]
\begin{center}
\leavevmode
\epsfxsize =2.7cm
\epsffile{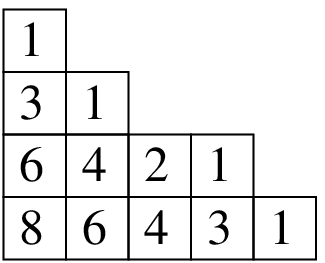}
\end{center}
\caption{\label{FIG1} The hook-lengths of $\alpha=(5,4,2,1)$}
\end{figure}

The criterion of irreducibility  
takes an especially nice form when $\alpha = \beta$. 
Let $\alpha'=(d_1,\ldots ,d_k)$ denote the conjugate partition.
Recall that the {\em hook-length} of the cell $c_{ij}$ of the Young diagram of
$\alpha$ lying on row $i$ and column $j$ is defined by
\[
h_{ij} = \alpha_i + d_j -i -j +1.
\]
This is illustrated in Figure~\ref{FIG1}.
\begin{corollary}\label{COR27}
The induction product $S_\alpha(t^a)\odot S_\alpha(t^b)$
is irreducible if and only if $|b-a|$ is not one of the hook-lengths
of $\alpha$.
\end{corollary}
\proof
Since $\alpha = \beta$, the sets $\I$ and $\J$ simply differ by a
translation of $b-a$.
Let $\I_0 := \I \setminus (\I \cap \J)$ and $\J_0 := \J \setminus (\I \cap \J)$.
If $b-a>0$ then $|\J_0| > |\I_0|$, $a+r+1,\,b+\alpha_1$ 
belong to $\J_0$ and $\I_0 \subset ]a+r+1,\,b+\alpha_1[$.
Therefore $\I$ and $\J$ are separated if and only if $\I_0=\emptyset$,
and in this case, they are strongly separated.
It is now a simple exercise to verify that $\I \subset \J$
if and only if $b-a$ is not one of the hook-lengths of $\alpha$.
\cqfd

\subsection{}
Our proof of Theorem~\ref{MAIN1} was based on the theory of canonical
bases for the quantum algebras $U_v(\Sl_N)$. 
There, $v$ was regarded as a formal variable.
To prove Theorem~\ref{MAIN2} we will employ the representation theory of the
quantum enveloping algebras $U_v(\slchap_N)$ of the Kac-Moody
Lie algebras $\slchap_N$. But in this proof we will set
$v=t^{\frac12}$.

\subsubsection
The link between the representation theories of the affine Hecke
algebras $\H_m=\H_m(t)$
and the quantum enveloping algebras $U_v(\slchap_N)$ 
was discovered by Drinfeld \cite{D}, in the degenerate case.
In the non-degenerate case this link was established by Cherednik in \cite{Ch}.
We will use the version of this link due to Chari and Pressley \cite{CP}.
In particular, we will take the definition of the Hopf algebra
$U_v(\slchap_N)$ from \cite{CP}, 2.1 and 2.4.
Note that $U_v(\slchap_N)$ contains $U_v(\Sl_N)$ as a Hopf subalgebra.
Our definition of the affine
Hecke algebra $\H_m$ coincides with the definition from \cite{CP}, 3.1. 
Fix the integer $N\geqslant2$. Set $\upsilon=t^{\frac12}$.

For each $m=1,2,\,\ldots$
there is a functor ${\cal D}$ from the category of all finite-dimensional 
$\H_m$-modules to the category of finite-dimensional 
$U_v(\slchap_N)$-modules which are of level $m$ as
$U_v(\Sl_N)$-modules
(\cite{CP}, 4.2). Recall that a $U_v(\Sl_N)$-module is said to be of
level $m$ if all its irreducible components occur in the $m$th tensor
power of the vector representation of $U_v(\Sl_N)$.)  
Moreover, under the assumption $N>m$, the functor ${\cal D}$ is 
an equivalence of categories.
The proofs of all these statements have been given in \cite{CP}.

Denote by $V_\alpha(z)$ the image under the functor ${\cal D}$
of the evaluation $\H_m$-module $S_\alpha(z)$. This is an irreducible module
over the quantum affine algebra $U_v(\slchap_N)$, which is also called
an evaluation module (\cite{CP}, 5.4). Consider also the image $V_\beta(w)$
of the $\H_n$-module $S_\beta(w)$. Suppose that $N>m,n$. Then the
modules $V_\alpha(z),V_\beta(w)$ are non-zero.
By \cite{CP}, 4.4 the $U_v(\slchap_N)$-module
${\cal D}(S_\alpha(z)\odot S_\beta(w))$ is equivalent to the tensor
product $V_\alpha(z)\otimes V_\beta(w)$.

By {\cite{Zel80}, 8.7 the induction products
$S_\alpha(z)\odot S_\beta(w)$ and $S_\beta(w)\odot S_\alpha(z)$
are irreducible and equivalent whenever the ratio $z/w$ does not belong to 
$t^{\Z}=v^{2{\Z}}$. Then the $U_v(\slchap_N)$-modules
$V_\alpha(z)\otimes V_\beta(w)$ and $V_\beta(w)\otimes V_\alpha(z)$
are also irreducible and equivalent. So then there exists
a non-zero $U_v(\slchap_n)$-intertwining operator
$$
I_{\alpha\beta}(z,w):
V_\alpha(z)\otimes V_\beta(w)
\longrightarrow
V_\beta(w)\otimes V_\alpha(z)\,,
$$
unique up to a multiplier from $\C^\ast$.
By definition, the vector spaces of the modules $V_\alpha(z)$
are the same for different values of the parameter $z$, and so are
the spaces of the modules $V_\beta(w)$ for different values of $w$.
The multipliers from $\C^\ast$ can be so chosen
that the linear operator $I_{\alpha\beta}(z,w)$ depends on $z,w$
as a rational function of $z/w$.
Such a choice can be made by using the explicit realizations of the
modules $V_\alpha(z)$ and $V_\beta(w)$ from [{\bf Ch}], Proposition~1.5.
Assume that such a choice has been made, and write
$I_{\alpha\beta}(z,w)=I_{\alpha\beta}(z/w)$.
In the 
physical literature the operator $I_{\alpha\beta}(z/w)$
is called the trigonometric $\check R$-matrix
corresponding to the evaluation modules 
$V_\alpha(z)$ and $V_\beta(w)$, see for instance \cite{KMT,MT}.

\subsubsection{}
When $N>m+n$, the irreducibility of the $U_v(\slchap_N)$-module
$V_\alpha(z)\otimes V_\beta(w)$ is equivalent to the irreducibility of the
$\H_{m+n}$-module $S_\alpha(z)\odot S_\beta(w)$. Thus for $N>m+n$ the
module $V_\alpha(z)\otimes V_\beta(w)$ is reducible, if and only if
$z/w=t^c$ for $c\in\Z$ described by our Theorem~\ref{MAIN1}. For an
arbitrary $N$, an irreducibility criterion for
$V_\alpha(z)\otimes V_\beta(w)$ has been recently given by Molev in
\cite{Molev}, independently of our results and by different methods;
see Theorems 3.1 and 4.1 therein.
Actually, Molev considers representations of the Yangian of $\gl_N$,
but it is known that the finite dimensional representation theory of 
the quantized affine algebra and that of the Yangian coincide,
see \cite{V}.  

\subsubsection{}
{\bf Proof of Theorem~\ref{MAIN2}.\ }
For each $k=1,\ldots,r$ consider the evaluation module
$S_{\alpha^{(k)}}(z_k)$. This is an irreducible module over
the affine Hecke algebra $\H_{m_k}$ where $\alpha^{(k)}$ is a partition of
$m_k$. Take the irreducible $U_v(\slchap_N)$-module
$V_k=V_{\alpha^{(k)}}(z_k)$. Suppose $N>m_k$ for every $k$, so that
$V_k$ is not zero. The
$U_v(\slchap_N)$-module
${\cal D}(S_{\alpha^{(1)}}(z_1)\odot\dots\odot S_{\alpha^{(r)}}(z_r))$
is equivalent to the tensor product
$$
{\cal D}(S_{\alpha^{(1)}}(z_1))
\otimes\dots\otimes
{\cal D}(S_{\alpha^{(r)}}(z_r))
=V_1\otimes\dots\otimes V_r\,.
$$
We will demonstrate that
the $U_v(\slchap_N)$-module $V_1\otimes\dots\otimes V_r$ is
irreducible,
if and only if the tensor products $V_k\otimes V_l$ are irreducible for
all $k<l$. When we choose $N>m_1+\dots+m_r$ here,
Theorem~\ref{MAIN2} will follow.  

Let $k,l=1,\dots,r$. By taking the first
non-zero coefficient of the Laurent expansion of 
the rational function $I_{\alpha^{(k)}\alpha^{(l)}}(u)$
at $u\to z_k/z_l$, we obtain an $U_v(\slchap_N)$-intertwining operator
$
I_{kl}:\ V_k\otimes V_l\longrightarrow V_l\otimes V_k.
$
By Corollary 38 the $\H_{2m_k}$-module
$
S_{\alpha^{(k)}}(z^{(k)})\odot S_{\alpha^{(k)}}(z^{(k)})
$ 
is irreducible for every $k=1,\dots,r$. So the tensor square
$V_k\otimes V_k$
is also irreducible. 
This implies that $I_{kk}$ equals the unit operator in
$V_k\otimes V_k$, up to a multiplier. A different proof
of this equality has been given in \cite{DJKM}.
Under these unit
conditions the argument of \cite{KMT,MT} shows, that
the irreducibility of the tensor
product $V_1\otimes\dots\otimes V_r$ is equivalent to the invertibility
of all operators
$I_{kl}$ with $k<l$. Here we again use the explicit realizations 
of the modules $V_k$ from 
[{\bf Ch}], Proposition~1.5. For the details of this argument
see {\cite{NT2}.
Applying this result in the particular case $r=2$,
we obtain that the invertibility of the single operator $I_{kl}$ with $k<l$
is equivalent to the irreducibility of the tensor product $V_k\otimes V_l$.
This completes the proof of Theorem~\ref{MAIN2}.

\subsubsection{}
We note that in the case where all the partitions $\alpha^{(k)}$
are rectangular, the irreducibility criterion for the module
$V_1\otimes\dots\otimes V_r$ was obtained in \cite{NT}. In particular,
it was shown in \cite{NT} that for the rectangular partitions
$\alpha^{(k)}$,
the module $V_1\otimes\dots\otimes V_r$ is irreducible, if and only if
the pairwise tensor products $V_k\otimes V_l$ are irreducible for all
$k<l$.
In fact, \cite{NT} discusses representations of the Yangian of $\gl_N$,
but as mentioned above, 
the finite dimensional representation theory of $U_v(\slchap_N)$
and that of the Yangian coincide.

\subsubsection{}
Let us make an important remark on the trigonometric
$\check R$-matrix $I_{\alpha\beta}(u)$. This is an operator-valued
rational function of $u$, determined by two partitions $\alpha$ and
$\beta$ of $m$ and $n$ respectively. It also depends on the choice of
$N>m,n$.
By definition, all zeroes and poles of this rational function of $u$
belong to $t^{\Z}$. Let us call a point $c\in\Z$ {\it singular,\/}
if the first non-zero coefficient in the Laurent expansion of 
$I_{\alpha\beta}(u)$ at $u\to t^c$ is a non-invertible linear operator.
Our proof of Theorem~\ref{MAIN2} shows, that for $N>m+n$
the point $c\in\Z$ is singular, if and only if the induced
module $S_\alpha(z)\odot S_\beta(w)$ with $z/w=t^c$ is reducible.
Thus when $N>m+n$, the list of all singular points $c$ for $I_{\alpha\beta}(u)$
is also given by Theorem~\ref{MAIN1}. 
In the particular case when every part of the partitions
$\alpha$ and $\beta$ is $1$, this result follows from \cite{AK}.

\subsection{}
{\bf Proof of Theorem~\ref{MainTH}.\ }
Consider the product 
\[
\pi=\langle{\cal I}_1\rangle\cdots\langle{\cal I}_r\rangle=
G^\ast(\m_1)\cdots G^\ast(\m_r)
\]
of $r$ quantum flag minors.
By Theorem~\ref{THAR} and Proposition~\ref{PROP15},
the normalized product $v^b\pi$ belongs for some $b\in\Z$
to the canonical basis ${\cal B}_v$ if and only if
the induction product $L=L_{\m_1}\odot\cdots\odot L_{\m_r}$
is irreducible. 
By Theorem~\ref{MAIN2}, the irreducibility of
$L$ is equivalent to the irreducibility of
all pairwise induction products $L_{\m_k}\odot L_{\m_l}$ where
$1\le k<l\le r$.

Again by Theorem~\ref{THAR} and Proposition~\ref{PROP15}, the induction product
$L_{\m_k}\odot L_{\m_l}$ is irreducible, if and only if
the product of two quantum flag minors
$G^\ast(\m_k)G^\ast(\m_l)=\langle{\cal I}_k\rangle\langle{\cal I}_l\rangle$
belongs to $\cal B_\upsilon$ up to a power of $v$.
But by Theorem~\ref{corol28}, this happens if and only the quantum flag minors
$\langle{\cal I}_k\rangle$ and $\langle{\cal I}_l\rangle$ quasi-commute.

\bigskip
{\bf Acknowledgements.}
We would like to express our gratitude to N. Kitanine for drawing our
attention to the results of \cite{KMT, MT}. Support from the EPSRC
and from the EC under the TMR grant FMRX-CT97-0100 is gratefully
acknowledged.


\bigskip
\small

\noindent
\begin{tabular}{ll}
{\sc B.  Leclerc} : &
D\'epartement de Math\'ematiques,
Universit\'e de Caen, Campus II,\\
& Bld Mar\'echal Juin,
BP 5186, 14032 Caen cedex, France\\
&email : {\tt leclerc@math.unicaen.fr}\\[5mm]
{\sc M. Nazarov} :&
Department of Mathematics,
University of York,\\
& York YO10 5DD, England\\
&email : {\tt mln1@york.ac.uk}\\[5mm]
{\sc J.-Y. Thibon} :&
Institut Gaspard Monge,
Universit\'e de Marne-la-Vall\'ee,\\
&Champs-sur-Marne,
77454 Marne-la-Vall\'ee cedex 2, France\\
&email : {\tt jyt@weyl.univ-mlv.fr}
\end{tabular}
\end{document}